\documentclass{amsart}
\usepackage{amscd,amssymb,enumerate,amsmath}
\usepackage[all]{xy}

\newcommand{\ppq}{\leqslant}
\newcommand{\pgq}{\geqslant}

\newcommand{\dsty}{\displaystyle}

\newcommand{\kk}{\mathbb{K}}
\newcommand{\zz}{\mathbb{Z}}

\newcommand{\id}{\mathrm{id}}

\newcommand{\im}{\operatorname{Im}\nolimits}
\renewcommand{\ker}{\operatorname{Ker}\nolimits}
\newcommand{\car}{\operatorname{char}\nolimits}
\newcommand{\Hom}{\operatorname{Hom}\nolimits}
\newcommand{\Ext}{\operatorname{Ext}\nolimits}
\newcommand{\HH}{\operatorname{HH}\nolimits}
\newcommand{\opp}{\operatorname{op}\nolimits}

\newcommand{\ev}{\operatorname{ev}\nolimits}

\newcommand{\ot}{\otimes}
\newcommand{\lan}{\Lambda_N}
\newcommand{\la}{\Lambda_1}
\renewcommand{\L}{\Lambda}
\newcommand{\oa}{\bar{a}}
\newcommand{\mo}{\mathfrak{o}}
\newcommand{\mt}{\mathfrak{t}}
\newcommand{\rrad}{\mathfrak{r}}

\newcommand{\N}{\mathcal{N}}

\newcommand{\md}[1]{\equiv{#1}\!\!\!\pmod{4}}
\newcommand{\set}[1]{\left\{ #1 \right\}}
\newcommand{\abs}[1]{\left|{#1}\right|}
\newcommand{\rep}[1]{\langle{#1}\rangle}

\newtheorem{theorem}{{Theorem}}[section]
\newenvironment{thm}{\begin{theorem}}{\end{theorem}}
\newcommand{\bt}{\begin{thm}}
\newcommand{\et}{\end{thm}}

\newtheorem{corollaire}[theorem]{{Corollary}}
\newenvironment{cor}{\begin{corollaire}}{\end{corollaire}}
\newcommand{\bc}{\begin{cor}}
\newcommand{\ec}{\end{cor}}

\newtheorem{lemme}[theorem]{{Lemma}}
\newenvironment{lemma}{\begin{lemme}}{\end{lemme}}
\newcommand{\bl}{\begin{lemma}}
\newcommand{\el}{\end{lemma}}

\newtheorem{proposition}[theorem]{{Proposition}}
\newcommand{\bprop}{\begin{proposition}}
\newcommand{\eprop}{\end{proposition}}

\newtheorem{definition}[theorem]{{Definition}}
\newenvironment{dfn}{\begin{definition} \rm}{\end{definition}}
\newcommand{\bd}{\begin{dfn}}
\newcommand{\ed}{\end{dfn}}

\newtheorem*{remark}{Remark}
\newcommand{\br}{\begin{remark}}
\newcommand{\er}{\end{remark}}

\begin{document}

\topmargin 0cm
\oddsidemargin 0.5cm
\evensidemargin 0.5cm
\pagestyle{plain}

\title[The Hochschild cohomology ring \dots]
{The Hochschild cohomology ring of a class of special biserial algebras}
\author[Snashall]{Nicole Snashall}
\address{Nicole Snashall\\Department of Mathematics\\
University of Leicester\\
University Road\\
Leicester, LE1 7RH\\
England}
\email{N.Snashall@mcs.le.ac.uk}
\author[Taillefer]{Rachel Taillefer}
\address{Rachel Taillefer\\Laboratoire de Math\'ematiques\\
Facult\'e des Sciences et Techniques\\
23 Rue Docteur Paul Michelon\\
42023 Saint-Etienne Cedex 2\\
France}
\email{rachel.taillefer@univ-st-etienne.fr}
\subjclass[2000]{Primary: 16E40, 18G10, 16D40, 16D50. Secondary:  16S37. }

\begin{abstract}
We consider a class of self-injective special biserial algebras
$\L_N$ over a field $K$ and show that the Hochschild cohomology ring
of $\L_N$ is a finitely generated $K$-algebra. Moreover the
Hochschild cohomology ring of $\L_N$ modulo nilpotence is a finitely
generated commutative $K$-algebra of Krull dimension two. As a
consequence the conjecture of \cite{SS}, concerning the Hochschild
cohomology ring modulo nilpotence, holds for this class of algebras.
\end{abstract}

\date{\today}
\maketitle

\section*{Introduction}

Let $\kk$ be a field. For $m \pgq 1$, let ${\mathcal Q}$ be the
quiver with $m$ vertices, labeled $0, 1 \ldots, m-1$, and $2m$
arrows as follows:
$$\xymatrix@=.01cm{
&&&&&&&&&&&\cdot\ar@/^.5pc/[rrrrrd]^{a}\ar@/^.5pc/[llllld]^{\bar{a}}\\
&&&&&&\cdot\ar@/^.5pc/[rrrrru]^{a}\ar@/^.5pc/[llldd]^{\bar{a}}&&&&&&&&&&\cdot\ar@/^.5pc/[rrrdd]^{a}\ar@/^.5pc/[lllllu]^{\bar{a}}\\\\
&&&\cdot\ar@/^.5pc/[rrruu]^{a}\ar@{.}@/_.3pc/[ldd]&&&&&&&&&&&&&&&&\cdot\ar@/^.5pc/[llluu]^{\bar{a}}\ar@{.}@/^.3pc/[rdd]\\\\
&&&&&&&&&&&&&&&&&&&&&&&\\
\\
\\
\\\\\\\\\\\\\\\\
&&&&&&&&&&&&&&&&\\
&&&&&&&&&&&\cdot\ar@/_.3pc/@{.}[rrrrru]\ar@/^.3pc/@{.}[lllllu] }$$
Let $a_i$ denote the arrow that goes from vertex $i$ to vertex
$i+1$, and let $\oa_i$ denote the arrow that goes from vertex $i+1$
to vertex $i$, for each $i=0, \ldots , m-1$ (with the obvious
conventions modulo $m$). We denote the trivial path at the vertex
$i$ by $e_i$. Paths are written from left to right.

In this paper, we study the Hochschild cohomology ring of a family
of algebras $\L_N$ given by this quiver ${\mathcal Q}$ and certain
relations. Let $N \pgq 1$, and let $\lan=\kk{\mathcal Q}/I_N$ where
$I_N$ is the ideal of $\kk{\mathcal Q}$ generated by
$$a_ia_{i+1},\ \ \oa_{i-1}\oa_{i-2},\ \
(a_i\oa_i)^{N}-(\oa_{i-1}a_{i-1})^{N}, \ \ \mbox{ for } i=0, 1,
\ldots , m-1,$$ where the subscripts are taken modulo $m$. These
algebras are all self-injective special biserial algebras and as
such play an important role in various aspects of representation
theory of algebras.

First, consider the case $N=1$, so that $\L_1 = \kk{\mathcal Q}/I_1$
where $I_1$ is the ideal of $\kk{\mathcal Q}$ generated by
$a_ia_{i+1}$, $\oa_{i-1}\oa_{i-2}$ and $a_i\oa_i-\oa_{i-1}a_{i-1}$,
for $i=0, \ldots , m-1$. In the case where $m$ is even, this algebra
occurred in the presentation by quiver and relations of the Drinfeld
double of the generalised Taft algebras studied in \cite{EGST}. It
also occurs in the study of the representation theory of
$U_q(\mathfrak{sl}_2)$; see \cite{Patra, Suter, Xiao}; see also
\cite{CK}. In the case $N=1$, the Hochschild cohomology ring of the
algebra $\L_1$ with $m=1$ is discussed in \cite{BGMS}, where
$q$-analogues of $\L_1$ were used to answer negatively a question of
Happel, in that they have finite dimensional Hochschild cohomology
ring but are of infinite global dimension when $q\in\kk^*$ is not a
root of unity. The more general algebras $\L_N$ occur in \cite{FS,
FS2}, in which the authors determine the Hopf algebras associated to
infinitesimal groups whose principal blocks are tame when $\kk$ is
an algebraically closed field of characteristic $p \pgq 3$; among
the principal blocks that are obtained in this classification are
the algebras $\L_{p^n}$ with $m=p^r$ vertices, for some integers $n$
and $r$.

This paper describes explicitly the structure of the Hochschild
cohomology ring of the algebras $\L_N$, for all $N \pgq 1$ and $m
\pgq 1$, and in all characteristics. In particular, we determine the
Hochschild cohomology ring of the algebras $\L_1$ of \cite{EGST} and
the algebras of \cite{FS}.

The main results show that the Hochschild cohomology ring,
$\HH^*(\L_N)$, is a finitely generated $\kk$-algebra.  We give an
explicit basis for each cohomology group $\HH^n(\L_N)$ together with
generators of the algebra $\HH^*(\L_N)$. We then determine
explicitly the Hochschild cohomology ring modulo nilpotence,
$\HH^*(\L_N)/\N$, where $\N$ is the ideal of $\HH^*(\L_N)$ generated
by all homogeneous nilpotent elements. Furthermore, we show that the
conjecture of \cite{SS} holds for all $\L_N$, that is, that
$\HH^*(\L_N)/\N$ is a finitely generated $\kk$-algebra. Moreover, we
show that $\HH^*(\L_N)/\N$ is a commutative ring of Krull dimension
2.

It is to be expected that the results in this paper will give some
insight into the more general problem of verifying the conjecture of
\cite{SS} for all special biserial algebras, and, indeed, for all
Koszul algebras since the algebra $\L_1$ is a Koszul algebra. A
general description of the multiplicative structure of the
Hochschild cohomology ring of a Koszul algebra was given in
\cite{BGSS}, but it remains open as to whether the conjecture of
\cite{SS} holds for all Koszul algebras. For comparison, we recall
that the conjecture has been verified for finite-dimensional
monomial algebras $A$ in \cite{GSS}, where it was shown that
$\HH^*(A)/\N$ is a finitely generated commutative algebra of Krull
dimension at most 1.

\bigskip

We now outline the structure of the paper. Section 1 describes the
minimal projective resolution $(P^*,\partial^*)$ of the algebra
$\L_N$ as a bimodule, together with an explanation of the
construction of this resolution. This construction is of more
general interest and applicability. Section 2 describes the methods
used to find the dimensions of the kernel and image of the induced
maps in the complex $\Hom(P^*,\L_N)$ and hence gives the dimension
of each Hochschild cohomology group $\HH^n(\L_N)$ in the case $m
\pgq 3$. In calculating $\HH^*(\L_N)$, we first consider the general
case $m \pgq 3$, and start with a description in section 3 of the
centre of $\L_N$, that is, of $\HH^0(\L_N)$, before giving the
structure of the Hochschild cohomology ring for $m \pgq 3$ in
sections 4 and 5. Section 6 deals with the case $m=2$ and section 7
with the case $m=1$. The final section summarises the paper, and remarks
that the finiteness conditions
(Fg1) and (Fg2) of \cite{EHSST} hold when $N=1$, thus enabling one to
describe the support varieties for finitely generated modules over the
algebras $\L_1$.

\bigskip

\section{A minimal projective bimodule resolution}

\bigskip

Let $N \pgq 1$, and let $\lan=\kk{\mathcal Q}/I_N$ where $I_N$ is
the ideal of $\kk{\mathcal Q}$ generated by $a_ia_{i+1},
\oa_{i-1}\oa_{i-2}$ and $(a_i\oa_i)^{N}-(\oa_{i-1}a_{i-1})^{N}$, for
$i =0, 1, \ldots , m-1$. Where there is no confusion, we label the
arrows of ${\mathcal Q}$ generically by $a$ and $\bar{a}$. We write
$\mo(\alpha)$ for the trivial path corresponding to the origin of
the arrow $\alpha$, so that $\mo(a_i) = e_i$ and $\mo(\oa_i) =
e_{i+1}$. We write $\mt(\alpha)$ for the trivial path corresponding
to the terminus of the arrow $\alpha$, so that $\mt(a_i) = e_{i+1}$
and $\mt(\oa_i) = e_i$. Recall that a non-zero element $r \in
\kk{\mathcal Q}$ is said to be uniform if there are vertices $v, w$
such that $r = vr = rw$. Let $S_i$ denote the simple module at the
vertex $i$, so that $\{S_0, \ldots , S_{m-1}\}$ is a complete set of
non-isomorphic simple modules for $\Lambda_N$.

The Hochschild cohomology ring of a $\kk$-algebra $\Lambda$ is given
by $\HH^*(\L) = \Ext^*_{\L^e}(\L,\L) = \oplus_{n \pgq
0}\Ext^n_{\L^e}(\L, \L)$ with the Yoneda product, where $\L^e =
\L^{\opp} \otimes_{\kk} \L$ is the enveloping algebra of $\L$. Since
all tensors are over the field $\kk$ we write $\otimes$ for
$\otimes_\kk$ throughout.

We now describe a minimal projective bimodule resolution $(P^*,
\partial^*)$ of $\Lambda_N$. By \cite{H}, we know that the multiplicity of
$\Lambda e_i\otimes e_j\Lambda$ as a direct summand of $P^n$ is
equal to the dimension of $\Ext^n_\Lambda(S_i, S_j)$. Thus we have,
for $n \pgq 0$, that
$$P^n = \oplus_{i=0}^{m-1} [\oplus_{r=0}^n \Lambda_N e_i\otimes
e_{i+n-2r}\Lambda_N].$$

We wish to label the summands of $P^n$ by certain elements in the
path algebra $\kk{\mathcal Q}$. The description given here is
motivated by \cite{GSn} where the first terms of a minimal bimodule
resolution of a finite-dimensional quotient of a path algebra were
determined explicitly from the early terms of the minimal right
$\Lambda$-resolution of $\Lambda/\rrad$, where $\rrad$ is the
Jacobson radical of $\L$, using \cite{GSZ}.

We start by recalling briefly the theory of projective resolutions
developed in \cite{GSZ} and \cite{GSn}. In \cite{GSZ}, the authors
give an explicit inductive construction of a minimal projective
resolution of $\Lambda/\rrad$ as a right $\Lambda$-module, for a
finite-dimensional algebra $\Lambda$ over a field $\kk$. For
$\Lambda = \kk{\mathcal Q}/I$ and finite-dimensional, they define
$g^0$ to be the set of vertices of ${\mathcal Q}$, $g^1$ to be the
set of arrows of ${\mathcal Q}$, and $g^2$ to be a minimal set of
uniform relations in the generating set of $I$, and then show that
there are subsets $g^n, n\ge 3$, of $\kk{\mathcal Q}$, where $x\in
g^n$ are uniform elements satisfying $x=\sum_{y\in g^{n-1}}yr_y
=\sum_{z\in g^{n-2}}zs_z$ for unique $r_y,s_z\in \kk{\mathcal Q}$,
which can be chosen in such a way that there is a minimal projective
$\Lambda$-resolution of the form
$$\cdots \to Q^4\to Q^3 \rightarrow Q^2 \rightarrow Q^1 \rightarrow Q^0
\rightarrow \Lambda/\rrad \rightarrow 0$$ having the following
properties.
\begin{itemize}
\item[(1)] For each $n\ge 0$,
$Q^n = \coprod_{x \in g^n} \mathfrak{t}(x)\Lambda$.
\item[(2)]
For each $x \in g^n$, there are unique elements $r_j \in
\kk{\mathcal Q}$ with $x = \sum_{j}g^{n-1}_jr_j$.
\item[(3)] For each $n\ge 1$,
using the decomposition of (2), for $x\in g^n$, the map $Q^n
\rightarrow Q^{n-1}$ is given by
$$\mathfrak{t}(x)\lambda \mapsto
\sum_{j}r_j\mathfrak{t}(x)\lambda .$$
\end{itemize}
where the elements of the set $g^n$ are labeled by $g^n =
\{g^n_j\}$. Thus the maps in this minimal projective resolution of
$\Lambda/\rrad$ as a right $\Lambda$-module are described by the
elements $r_j$ which are uniquely determined by (2).

In \cite{GSn}, these same sets $g^n$ are used to give an explicit
description of the first three maps in a minimal projective
resolution of $\Lambda$ as a $\Lambda,\Lambda$-bimodule, thus
connecting the minimal projective resolution of $\Lambda/\rrad$ as a
right $\Lambda$-module with a minimal projective resolution of
$\Lambda$ as a $\Lambda,\Lambda$-bimodule. We use the same ideas
here to give a minimal projective bimodule resolution for the
algebra $\Lambda_N$, for all $N \pgq 1$. In the case where $N=1$,
the algebra $\L_1$ is Koszul, and we refer to \cite{GHMS} which uses
this approach and gives a minimal projective bimodule resolution for
any Koszul algebra.

\bigskip

We start with the case $N=1$ and define sets $g^n$ in the path
algebra $\kk{\mathcal Q}$ which we will use to label the generators
of $P^n$.

\subsection{The bimodule resolution for $\Lambda_1$}

Consider the algebra $\L_1$. For each vertex $i = 0, \ldots , m-1$
and for each $r = 0, \ldots , n$, we define elements $g^n_{r,i}$ in
$\kk{\mathcal Q}$ as follows. Let
$$g^n_{r,i} = \sum_{p}(-1)^sp$$ where the sum is over all paths $p$
of length $n$, written $p = \alpha_1 \alpha_2 \cdots \alpha_n$ where
the $\alpha_i$ are arrows in ${\mathcal Q}$, such that
\begin{enumerate}
\item[(i)] $p \in e_i\kk{\mathcal Q}$,
\item[(ii)] $p$ contains $r$ arrows of the form $\bar{a}$ and $n-r$
arrows of the form $a$, and
\item[(iii)] $s = \sum_{\alpha_j = \bar{a}}j$.
\end{enumerate}
It follows that $g_{r,i}^n \in e_i(\kk{\mathcal Q})e_{i+n-2r}$, for
$i=0,\ldots, m-1$ and $r=0, \ldots, n$. Since the elements
$g^n_{r,i}$ are uniform elements, we may define
$\mathfrak{o}(g^n_{r,i}) = e_i$ and $\mathfrak{t}(g^n_{r,i}) =
e_{i+n-2r}$. Then $$P^n =
\oplus_{i=0}^{m-1}[\oplus_{r=0}^{n}\Lambda_1\mathfrak{o}(g^n_{r,i})
\otimes \mathfrak{t}(g^n_{r,i})\Lambda_1].$$

We set
$$g^n = \bigcup_{i=0}^{m-1}\{g^n_{r,i} \mid r = 0, \ldots , n\}.$$
It is easy to see, for the cases $n = 0, 1, 2$, that we have
$g_{0,i}^0 = e_i$, $g_{0,i}^1 = a_i$ and $g_{1,i}^1 =
-\bar{a}_{i-1}$, whilst $g_{0,i}^2 = a_ia_{i+1}$, $g_{1,i}^2 =
a_i\bar{a}_i - \bar{a}_{i-1}a_{i-1}$ and $g_{2,i}^2 =
-\bar{a}_{i-1}\bar{a}_{i-2}$. Thus
$$\begin{array}{lll}
g^0 & = & \{e_i \mid i = 0, \ldots , m-1\},\\
g^1 & = & \{a_i, -\oa_i \mid i = 0, \ldots , m-1\},\\
g^2 & = & \{a_ia_{i+1},\ \ a_i\oa_i-\oa_{i-1}a_{i-1},\ \ -\oa_{i-1}\oa_{i-2} \mbox{ for all $i$}\},\\
\end{array}$$
so that $g^2$ is a minimal set of uniform relations in the
generating set of $I_1$.

\bigskip

To describe the map $\partial^n\colon P^n\to P^{n-1}$, we need to
write the elements of the set $g^n$ in terms of the elements of the
set $g^{n-1}$. The proof of the next lemma is straightforward, and
is left to the reader.

\bigskip

\bl\label{lemma:maps} For the algebra $\L_1$, for $i =0, 1, \ldots ,
m-1$ and $r = 0, 1, \ldots , n$, we have:
$$\begin{array}{lll}
g^n_{r,i} & = g^{n-1}_{r,i}a_{i+n-2r-1} +
(-1)^ng^{n-1}_{r-1,i}\bar{a}_{i+n-2r} & = (-1)^ra_ig^{n-1}_{r,i+1} +
(-1)^r\bar{a}_{i-1}g^{n-1}_{r-1,i-1}\end{array}$$
with the
conventions that $g^n_{-1,i} = 0$ and $g^{n-1}_{n,i} = 0$ for all
$n, i$. Thus
$$g^n_{0,i} = g^{n-1}_{0,i}a_{i+n-1} = a_ig^{n-1}_{0,i+1} \mbox{ and }
g^n_{n,i} = (-1)^ng^{n-1}_{n-1,i}\bar{a}_{i-n+1} = (-1)^n\bar{a}_{i-1}g^{n-1}_{n-1,i-1}.$$
\el

\bigskip

Since $\L_1$ is Koszul, we now use \cite{GHMS} to give a minimal
projective bimodule resolution $(P^n, \partial^n)$ of $\L_1$. We
define the map $\partial^0\colon P^0\to \L_1$ to be the
multiplication map. However, to define $\partial^n$ for $n \pgq 1$,
we first need to introduce some notation. In describing the image of
$\mathfrak{o}(g^n_{r,i}) \otimes \mathfrak{t}(g^n_{r,i})$ under
$\partial^n$ in the projective module $P^{n-1}$, we use subscripts
under $\otimes$ to indicate the appropriate summands of the
projective module $P^{n-1}$. Specifically, let $\otimes_r$ denote a
term in the summand of $P^{n-1}$ corresponding to $g^{n-1}_{r,-}$,
and $\otimes_{r-1}$ denote a term in the summand of $P^{n-1}$
corresponding to $g^{n-1}_{r-1,-}$, where the appropriate index $-$
of the vertex may always be uniquely determined from the context.
Indeed, since the relations are uniform along the quiver, we can
also take labeling elements defined by a formula independent of $i$,
and hence we omit the index $i$ when it is clear from the context.
Recall that nonetheless all tensors are over $\kk$.

Now, for $n\pgq 1$, keeping the above notation and using
\cite{GHMS}, we define the map $\partial^n: P^n\to P^{n-1}$ for the
algebra $\Lambda_1$ as follows:
$$\begin{array}{ll}
\partial^n \colon\mathfrak{o}(g^n_{r,i}) \otimes \mathfrak{t}(g^n_{r,i})
\mapsto & (e_i\otimes_ra_{i+n-2r-1} +
(-1)^ne_i\otimes_{r-1}\bar{a}_{i+n-2r})\\
& +(-1)^n((-1)^ra_i\otimes_re_{i+n-2r} +
(-1)^r\bar{a}_{i-1}\otimes_{r-1}e_{i+n-2r}).
\end{array}$$
Using our conventions that $g^n_{-1,i} = 0$ and $g^{n-1}_{n,i} =
0$ for all $n, i$, the degenerate cases ($r = 0, r = n$) simplify to
$$\partial^n\colon \mathfrak{o}(g^n_{0,i}) \otimes \mathfrak{t}(g^n_{0,i})
\mapsto e_i\otimes_0a_{i+n-1} + (-1)^na_i\otimes_0e_{i+n}$$ where
the first term is in the summand corresponding to $g^{n-1}_{0,i}$
and the second term is in the summand corresponding to
$g^{n-1}_{0,i+1}$, whilst
$$\partial^n\colon \mathfrak{o}(g^n_{n,i}) \otimes \mathfrak{t}(g^n_{n,i})
\mapsto (-1)^ne_i\otimes_{n-1}\bar{a}_{i-n} +
\bar{a}_{i-1}\otimes_{n-1}e_{i-n},$$ with the first term in the
summand corresponding to $g^{n-1}_{n-1,i}$ and the second term in
the summand corresponding to $g^{n-1}_{n-1,i-1}$.

The following result is now immediate from \cite[Theorem 2.1]{GHMS}.

\bigskip

\bt
With the above notation, $(P^n,
\partial^n)$ is a minimal projective resolution of $\Lambda_1$ as a
$\Lambda_1,\Lambda_1$-bimodule. \et

\subsection{The bimodule resolution for $\Lambda_N$ with $N \pgq 1$}

We now consider the general case $N \pgq 1$. In this case we use the
approach of \cite{GHMS,GSn} to construct a minimal projective
resolution for $\L_N$. We keep the conventions that $g^n_{-1,i} = 0$
and $g^{n-1}_{n,i} = 0$, for all $n, i$, throughout the paper.

\bigskip

\bd For the algebra $\L_N$ ($N \pgq 1$), for $i = 0, \ldots , m-1$
and $r = 0, \ldots , n$, we define $g_{0,i}^0 = e_i$, and,
inductively for $n \pgq 1$,
$$g^n_{r,i} = \left \{
\begin{array}{ll}
g^{n-1}_{r,i}a + (-1)^ng^{n-1}_{r-1,i}\bar{a}(a\oa)^{N-1} & \mbox{ if } n-2r > 0;\\
&\\
g_{r,i}^{n-1}a(\oa a)^{N-1}+(-1)^ng_{r-1,i}^{n-1}\oa & \mbox{ if } n-2r < 0;\\
&\\
g_{r,i}^{n-1}a(\oa a)^{N-1} + g_{r-1,i}^{n-1}\oa(a\oa)^{N-1} &
\mbox{ if } n=2r.
\end{array}\right.$$
\ed

\br
\begin{enumerate}
\item If $N = 1$, the above definition reduces to that given for $\L_1$.
\item We have $g_{0,i}^1 = a_i$ and $g_{1,i}^1 = -\bar{a}_{i-1}$, for all
$i$. Also $g_{0,i}^2 = a_ia_{i+1}$, $g_{1,i}^2 = (a_i\bar{a}_i)^N -
(\bar{a}_{i-1}a_{i-1})^N$ and $g_{2,i}^2 =
-\bar{a}_{i-1}\bar{a}_{i-2}$, for all $i$. Hence $g^2$ is a minimal
set of uniform relations in the generating set of $I_N$.
\item The projectives in a minimal projective bimodule resolution of $\Lambda_N$ are
given by $$P^n =
\oplus_{i=0}^{m-1}[\oplus_{r=0}^n\L_N\mo(g^n_{r,i})\otimes\mt(g^n_{r,i})\L_N].$$
\end{enumerate}
\er

The next result is an analogue of Lemma \ref{lemma:maps} for the
cases $n-2r > 0$ and $n-2r < 0$. The proof of the Lemma is easy to
verify and is omitted.

\bigskip

\bl\label{lem:g} For the algebra $\L_N$ ($N\pgq 1$), for $i = 0, \ldots , m-1$,
$n \pgq 1$, and $r = 0, \ldots , n$, we have:
$$g^n_{r,i} = \left \{
\begin{array}{ll}
g^{n-1}_{r}a + (-1)^ng^{n-1}_{r-1}\bar{a}(a\oa)^{N-1} & =
(-1)^rag^{n-1}_{r} + (-1)^r\bar{a}(a\oa)^{N-1}g^{n-1}_{r-1} \mbox{\ \ if } n-2r > 0;\\
&\\
g_r^{n-1}a(\oa a)^{N-1} + (-1)^ng_{r-1}^{n-1}\oa & = (-1)^ra(\oa
a)^{N-1}g_r^{n-1}+(-1)^{r}\oa g_{r-1}^{n-1} \mbox{\ \ if } n-2r < 0.
\end{array} \right.
$$
\el

We now define the maps $\partial^n$.

\bigskip

\bd\label{maps} For the algebra $\L_N$ ($N \pgq 1$), and for $n \pgq
0$, we define maps $\partial^n$ as follows. For $n = 0$, the map
$\partial^0\colon P^0\to \L_N$ is the multiplication map. For $n
\pgq 1$, and $i = 0, 1, \ldots , m-1$, the map $\partial^n\colon
P^n\rightarrow P^{n-1}$ is given
by $\partial^n \colon e_i\ot_r e_{i+n-2r}\mapsto$\\
$$\left \{
\begin{array}{ll}
e_i\ot_r e_{i+(n-1)-2r}a + (-1)^{n+r}ae_{i+1}\ot_r e_{i+1+(n-1)-2r} \\
\hspace{1cm} + (-1)^{n+r}\oa(a\oa)^{N-1}e_{i-1}\ot_{r-1} e_{i-1+(n-1)-2(r-1)}\\
\hspace{2cm} + (-1)^ne_i\ot_{r-1} e_{i+(n-1)-2(r-1)}\oa(a\oa)^{N-1} &\mbox{\ \ if $n-2r>0$,}\\
\ \\
e_i\ot_r e_{i+(n-1)-2r}a(\oa a)^{N-1} + (-1)^{n+r}a(\oa a)^{N-1}e_{i+1}\ot_r e_{i+1+(n-1)-2r} \\
\hspace{1cm} + (-1)^{n+r}\oa e_{i-1}\ot_{r-1} e_{i-1+(n-1)-2(r-1)}\\
\hspace{2cm} + (-1)^ne_i\ot_{r-1} e_{i+(n-1)-2(r-1)}\oa
&\mbox{\ \ if $n-2r<0$,}\\
\ \\
\sum_{k=0}^{N-1}(\oa a)^k[e_i\ot_{\frac{n}{2}} e_{i-1}a +
(-1)^{\frac{n}{2}}\oa e_{i-1}\ot_{\frac{n-2}{2}} e_i](\oa a)^{N-k-1} \\
\hspace{1cm} + (a\oa)^k[(-1)^{\frac{n}{2}}ae_{i+1}\ot_{\frac{n}{2}}
e_i + e_i\ot_{\frac{n-2}{2}} e_{i+1}\oa](a\oa)^{N-k-1} &\mbox{\ \ if
$n-2r=0$}.
\end{array}\right.$$
\ed

In order to prove that $(P^n, \partial^n)$ is a minimal projective
resolution of $\Lambda_N$ as a $\L_N,\L_N$-bimodule, we use an argument
which was given in \cite{GSn}. To do this, we note that
$$\L_N/\rrad\otimes_{\L_N}P^n \cong
\oplus_{i=0}^{m-1}\oplus_{r=0}^n\mt(g^n_{r,i})\L_N$$
as right $\L_N$-modules and that the map
$\id\otimes\partial^n : \L_N/\rrad\otimes_{\L_N}P^n \to
\L_N/\rrad\otimes_{\L_N}P^{n-1}$ is equivalent to the map
$\oplus_{i=0}^{m-1}\oplus_{r=0}^n\mt(g^n_{r,i})\L_N \to
\oplus_{i=0}^{m-1}\oplus_{r=0}^{n-1}\mt(g^{n-1}_{r,i})\L_N$ given by
$$\mt(g^n_{r,i})\mapsto
\left \{
\begin{array}{ll}
\mt(g^{n-1}_{r,i})a + (-1)^n\mt(g^{n-1}_{r-1,i})\oa(a\oa)^{N-1} &
\mbox{\ \ if $n-2r>0$,}\\
\mt(g^{n-1}_{r,i})a(\oa a)^{N-1} + (-1)^n\mt(g^{n-1}_{r-1,i})\oa &
\mbox{\ \ if $n-2r<0$,}\\
\mt(g^{n-1}_{\frac{n}{2},i})a(\oa a)^{N-1} +
\mt(g^{n-1}_{\frac{n-2}{2},i})\oa(a\oa)^{N-1} &
\mbox{\ \ if $n-2r=0$}.
\end{array}\right.$$
It is then straightforward (with Lemma \ref{lem:g}) to see that
$(\L_N/\rrad\otimes_{\L_N}P^n, \id\otimes\partial^n)$ is a minimal
projective resolution of $\Lambda_N/\rrad$ as a right
$\Lambda_N$-module, which satisfies the conditions of \cite{GSZ} as
explained above, for all $N \pgq 1$.
For completeness in the proof of Theorem \ref{thm:pr},
we explain in full the strategy used in \cite[Proposition 2.8]{GSn}.

\bigskip

\bt\label{thm:pr} With the above notation, $(P^n, \partial^n)$ is a minimal
projective resolution of $\Lambda_N$ as a
$\Lambda_N,\Lambda_N$-bimodule, for all $N \pgq 1$. \et

{\it Proof }  The first step is to verify that we have a complex by
showing that $\partial^2 = 0$; this is straightforward and the
details are left to the reader.

Now suppose for contradiction that $\ker\partial^{n-1}\not\subseteq
\im\partial^n$ for some $n \pgq 1$. Then there is some non-zero map
$\ker\partial^{n-1} \to \ker\partial^{n-1}/\im\partial^n$. Hence there
is a simple $\L_N,\L_N$-bimodule $U\otimes V$ and epimorphism
$\ker\partial^{n-1}/\im\partial^n \twoheadrightarrow U \otimes V$
such that the composition
$f: \ker\partial^{n-1} \to \ker\partial^{n-1}/\im\partial^n
\twoheadrightarrow U\otimes V$
is a non-zero epimorphism.

Now, $U$ is a simple left $\L_N$-module, $V$ is a simple right $\L_N$-module
and the functor $\L_N/\rrad \otimes_{\L_N} -$ preserves monos and epis.
Since $(\L_N/\rrad\otimes_{\L_N}P^n, \id\otimes\partial^n)$ is a minimal
projective resolution of $\Lambda_N/\rrad$ as a right
$\Lambda_N$-module, we have that
$\L_N/\rrad\otimes_{\L_N}\im\partial^n \cong \im(\id\otimes\partial^n) =
\ker(\id\otimes\partial^{n-1}) \cong
\L_N/\rrad\otimes_{\L_N}\ker\partial^{n-1}$.
Thus, applying the functor $\L_N/\rrad\otimes_{\L_N}-$ to the map
$\partial^n :P^n \to P^{n-1}$ gives:
$$\xymatrix{
\L_N/\rrad \otimes_{\L_N}
P^n\ar[d]^{\id\otimes\partial^n}\ar[r]^{\!\!\!\!\id\otimes\partial^n} &
\L_N/\rrad\otimes_{\L_N}P^{n-1} & \\
\L_N/\rrad\otimes_{\L_N}\im\partial^n\ar[r]^{\hspace*{-.5cm}\cong} &
\L_N/\rrad\otimes_{\L_N}\ker\partial^{n-1}\ar[r]^{\id \otimes f} &
\L_N/\rrad\otimes_{\L_N}U\otimes V.
}$$
The map $\id\otimes f$ is non-zero so the composition
$\L_N/\rrad\otimes_{\L_N}P^n \to \L_N/\rrad\otimes_{\L_N}U\otimes V$ is
non-zero. However, this is simply the functor
$\L_N/\rrad\otimes_{\L_N}-$
applied to the composition $P^n \stackrel{\partial^n}{\to}
\ker\partial^n \stackrel{f}{\to} U\otimes V$. Now the composition
$P^n \stackrel{\partial^n}{\to}
\ker\partial^n \stackrel{f}{\to} U\otimes V$
is zero since $\im\partial^{n-1}\subseteq\ker\partial^n$. This gives the required
contradiction, so the complex is exact.

Finally, minimality follows since we know
that the projectives are those of a minimal projective resolution of
$\L_N$ as a $\L_N,\L_N$-bimodule
from \cite{H}. \qed

\bigskip

We note that the degenerate cases $(r = 0, r = n)$ in the maps
$\partial^n$ may be written as:
$$\partial^n \colon e_i\ot_0 e_{i+n}\mapsto
e_i\ot_0 e_{i+n-1}a + (-1)^nae_{i+1}\ot_0 e_{i+n}$$ and
$$\partial^n \colon e_i\ot_n e_{i-n}\mapsto
\oa e_{i-1}\ot_{n-1} e_{i-n} + (-1)^ne_i\ot_{n-1} e_{i-n+1}\oa.
$$

\subsection{Conventions on notation}

So far, we have tried to simplify notation by denoting the
idempotent $\mathfrak{o}(g^n_{r,i}) \otimes \mathfrak{t}(g^n_{r,i})$
of the summand $\L_N\mathfrak{o}(g^n_{r,i}) \otimes
\mathfrak{t}(g^n_{r,i})\L_N$ of $P^n$ uniquely by $e_i\ot_r
e_{i+n-2r}$ where $0 \ppq i \ppq m-1$. However, even this notation
with subscripts under the tensor product symbol becomes cumbersome
in computations and in describing the elements of the Hochschild
cohomology ring. Thus we make some conventions (by abuse of
notation) which we keep throughout the paper.

First, we note that, in the language of \cite{BK}, the $n$-th
projective module $P^n$ in a projective bimodule resolution of
$\Lambda_N$ is denoted $S\otimes_SV_n\otimes_SS$ where $S$ is the
semisimple $\kk$-algebra generated by the idempotents $\{e_0, e_1,
\ldots , e_{m-1}\}$ and $V_n$ is the $S,S$-bimodule generated by the
set $$\{g^n_{r,i} \mid i=0, 1, \ldots , m-1,\ r = 0, 1, \ldots ,
n\}.$$ (Note that here $\otimes_S$ denotes a tensor over $S$, that
is, over a finite sum of copies of $\kk$; we continue to reserve the
notation $\otimes$ exclusively for tensors over $\kk$.) There is a
bijective correspondence between the idempotent
$\mathfrak{o}(g^n_{r,i}) \otimes \mathfrak{t}(g^n_{r,i})$ and the
term $e_i\otimes_Sg^n_{r,i}\otimes_S e_{i+n-2r}$ in $V_n$, for $i=0,
1, \ldots , m-1$ and $r = 0, 1, \ldots , n$.

Since $e_{i+n-2r} \in \{e_0, e_1, \ldots , e_{m-1}\}$, it would be
usual to reduce the subscript $i+n-2r$ modulo $m$. However, to make
it explicitly clear to which summand of the projective module $P^n$
we are referring and thus to avoid confusion, whenever we write
$e_i\otimes e_{i+k}$ for an element of $P^n$, we will always have $i
\in \{0, 1, \ldots , m-1\}$ and consider $i+k$ as an element of
${\mathbb Z}$, in that $r = (n-k)/2$ and $e_i\otimes e_{i+k} =
e_i\otimes_{\frac{n-k}{2}} e_{i+k}$ and thus lies in the
$\frac{n-k}{2}$-th summand of $P^n$. We do not reduce $i+k$ modulo
$m$ in any of our computations. In this way, when considering
elements in $P^n$, our element $e_i\otimes e_{i+k}$ corresponds
uniquely to the idempotent $\mathfrak{o}(g^n_{r,i}) \otimes
\mathfrak{t}(g^n_{r,i})$ of $P^n$ (or, equivalently, to
$e_i\otimes_Sg^n_{r,i}\otimes_S e_{i+n-2r}$) with $r = (n-k)/2$, for
each $i=0, 1, \ldots , m-1$.

Since we take $0\leq i \leq m-1$ we clarify our conventions in
Definition \ref{maps} in the cases $i = 0$ and $i = m-1$, to ensure
that we always have $0\leq j \leq m-1$ in the first idempotent entry
of each tensor $(e_j\ot -)$. Specifically, we have:\\
$\partial^n \colon e_0\ot e_{n-2r}\mapsto$\\
$$\left \{
\begin{array}{ll}
e_0\ot e_{(n-1)-2r}a + (-1)^{n+r}ae_{1}\ot e_{1+(n-1)-2r} \\
\hspace{1cm} + (-1)^{n+r}\oa(a\oa)^{N-1}e_{m-1}\ot
e_{m-1+(n-1)-2(r-1)}\\
\hspace{2cm} + (-1)^ne_0\ot e_{(n-1)-2(r-1)}\oa(a\oa)^{N-1} &\mbox{\ \ if $n-2r>0$,}\\
\ \\
e_0\ot e_{(n-1)-2r}a(\oa a)^{N-1} + (-1)^{n+r}a(\oa a)^{N-1}e_{1}\ot e_{1+(n-1)-2r} \\
\hspace{1cm} + (-1)^{n+r}\oa e_{m-1}\ot e_{m-1+(n-1)-2(r-1)}\\
\hspace{2cm} + (-1)^ne_0\ot e_{(n-1)-2(r-1)}\oa
&\mbox{\ \ if $n-2r<0$,}\\
\ \\
\sum_{k=0}^{N-1}(\oa a)^k[e_0\ot e_{-1}a +
(-1)^{\frac{n}{2}}\oa e_{m-1}\ot e_m](\oa a)^{N-k-1} \\
\hspace{1cm} + (a\oa)^k[(-1)^{\frac{n}{2}}ae_{1}\ot e_0 + e_0\ot
e_{1}\oa](a\oa)^{N-k-1} &\mbox{\ \ if $n-2r=0$,}
\end{array}\right.$$
and
$\partial^n \colon e_{m-1}\ot e_{m-1+n-2r}\mapsto$\\
$$\left \{
\begin{array}{ll}
e_{m-1}\ot e_{m-1+(n-1)-2r}a + (-1)^{n+r}ae_{0}\ot e_{(n-1)-2r} \\
\hspace{1cm} + (-1)^{n+r}\oa(a\oa)^{N-1}e_{m-2}\ot
e_{m-2+(n-1)-2(r-1)}\\
\hspace{2cm} + (-1)^ne_{m-1}\ot e_{m-1+(n-1)-2(r-1)}\oa(a\oa)^{N-1} &\mbox{\ \ if $n-2r>0$,}\\
\ \\
e_{m-1}\ot e_{m-1+(n-1)-2r}a(\oa a)^{N-1} + (-1)^{n+r}a(\oa a)^{N-1}e_0\ot e_{(n-1)-2r} \\
\hspace{1cm} + (-1)^{n+r}\oa e_{m-2}\ot e_{m-2+(n-1)-2(r-1)}\\
\hspace{2cm} + (-1)^ne_{m-1}\ot e_{m-1+(n-1)-2(r-1)}\oa &\mbox{\ \ if $n-2r<0$,}\\
\ \\
\sum_{k=0}^{N-1}(\oa a)^k[e_{m-1}\ot e_{m-2}a +
(-1)^{\frac{n}{2}}\oa e_{m-2}\ot e_{m-1}](\oa a)^{N-k-1} \\
\hspace{1cm} + (a\oa)^k[(-1)^{\frac{n}{2}}ae_{0}\ot e_{-1} +
e_{m-1}\ot e_m\oa](a\oa)^{N-k-1} &\mbox{\ \ if $n-2r=0$}.
\end{array}\right.$$

\bigskip

Our notation $e_i\otimes e_{i+k}$ is used without further comment
throughout the rest of the paper.

\bigskip

We end this section by determining the dimension of each space
$\Hom_{\Lambda_N^e}(P^n, \Lambda_N)$ for each $n \pgq 0$.

\bigskip

\subsection{The complex ${\Hom}_{\Lambda_N^e}(P^n, \Lambda_N)$ }

All our homomorphisms are $\Lambda_N^e$-homomorphisms and so we
write ${\Hom}(-,-)$ for ${\Hom}_{\Lambda_N^e}(-,-)$.

\medskip

\bl If $m \pgq 3$, write $n=pm+t$ where $p\pgq 0$ and $0\ppq t\ppq
m-1$. Then
$$\dim_\kk\Hom(P^n, \Lambda_N) =
\left\{
\begin{array}{ll}
(4p+2)mN & \mbox{if } t\neq m-1 \\
(4p+4)mN & \mbox{if } t=m-1.
\end{array}
\right.
$$
If $m = 1$ or $m = 2$ then $$\dim_\kk\Hom(P^n, \L_N) = 4N(n+1).$$
\el

{\it Proof } An element $f$ in ${\Hom}(P^n, \Lambda_N)$ is
determined by the elements $f(e_i\ot e_{i+n-2r})$, which can be
arbitrary elements of $e_i\Lambda_N e_{i+n-2r}$.

Suppose first that $m \pgq 3$. Then, for vertices $i, j$, we have
$$\dim e_i\Lambda_N e_j =
\left\{
\begin{array}{ll}
2N & \mbox{if } i = j,\\
N & \mbox{if } i-j\equiv \pm 1 \pmod{m},\\
0 & \mbox{otherwise.}
\end{array}
\right. $$ For a finite subset $X \subset \mathbb{Z}$ and $0 \ppq s
\ppq m-1$, we define $\nu_s(X) = |\{ x\in X: x\equiv s \pmod{m}\}|$.
With this we have
$$\dim {\Hom}(P^n, \Lambda_N) =
mN(2\nu_0(I_n) + \nu_1(I_n) + \nu_{-1}(I_n))
$$
where $I_n = \{ n-2r: 0\ppq r\ppq n\}$. We observe that the disjoint
union $I_n \cup I_{n-1}$ is equal to  $[-n,n]$, the set of all
integers $x$ with $-n\ppq x\ppq n$. Clearly $\nu_0([-n,n]) = 2p+1$
and
$$\nu_{1}([-n,n]) =
\left \{
\begin{array}{ll}
2p & \mbox{if } t=0 \\
2p+1 & \mbox{if } 1\ppq t \ppq m-2\\
2p+2 & \mbox{if } t=m-1.
\end{array} \right.
$$
Moreover $\nu_{-1}([-n,n]) = \nu_1([-n,n])$, by symmetry.

The lemma now follows by induction on $n$, with the case $n=0$ being
clear. For the inductive step, note that
$$\begin{array}{rl}
\dim{\Hom}(P^n, \Lambda_N) = & mN(2\nu_0(I_n) + \nu_1(I_n) +
\nu_{-1}(I_n)) \\
 = & mN(2\nu_0([-n,n]) + 2\nu_1([-n,n])) - \dim{\Hom}(P^{n-1},
 \Lambda_N).
\end{array}
$$
The statement for $\dim\Hom(P^n, \Lambda_N)$ now follows for $m \pgq
3$.

In the case $m = 1$, we have $P^n = \oplus_{r=0}^n\L_N\otimes \L_N$ and
$\dim\L_N = 4N$. Thus $\dim\Hom(P^n, \Lambda_N) = 4N(n+1)$.

Finally, for $m=2$, suppose first that $n$ is even. Then $P^n$ has
$n+1$ summands equal to $\L_N e_0\otimes e_0\L_N$ and $n+1$ summands
equal to $\L_N e_1\otimes e_1\L_N$. By symmetry, $\dim\Hom(P^n,
\Lambda_N) = 2(n+1)\dim e_1\L_N e_1 = 4N(n+1)$. For $n$ odd, $P^n$
has $n+1$ summands equal to $\L_N e_0\otimes e_1\L_N$ and $n+1$
summands equal to $\L_N e_1\otimes e_0\L_N$. Hence, again using
symmetry, we have $\dim\Hom(P^n, \Lambda_N) = 2(n+1)\dim e_0\L_N e_1
= 4N(n+1)$. This completes the proof. \qed

\bigskip

\section{The dimensions of the Hochschild cohomology spaces for $m\pgq3$.}

\bigskip

We give here the main arguments we used to compute the kernel and
the image of the differentials, and leave the details to the reader.
We denote both the map $P^n \to P^{n-1}$ and its induced map
$\Hom(P^{n-1},\L_N)\to\Hom(P^n,\L_N)$ by $\partial^n$. Then we have
that $\HH^n(\lan)=\ker(\partial^{n+1})/\im(\partial^n)$. We assume
without further comment that $m \pgq 3$ throughout this section.

\subsection{Explicit description of maps}

In order to compute the Hochschild cohomology groups, we need an
explicit description of the image of each $e_i\otimes_r e_{i+n-2r} =
e_i\otimes e_{i+n-2r}$ under $f \in \Hom(P^n, \Lambda_N)$ for $i =
0, 1, \ldots , m-1$. We write $n = pm + t$ with $0 \ppq t \ppq m-1$,
$p \pgq 0$.

The image $f(e_i\otimes e_{i+n-2r})$ lies in $e_i\L_N e_j$ where $j
\in {\mathbb Z}_m$ with $j \equiv i+n-2r \pmod{m}$. Since $e_i\L_N
e_j = (0)$ if the vertices $i$ and $j$ are sufficiently far apart in
the quiver of $\L_N$ (at least for $m \pgq 3$), for a non-zero image
$f(e_i\otimes e_{i+n-2r})$ we are only interested in examining the
cases where $i+n-2r \equiv i, i-1, i+1 \pmod{m}$. This leads to the
consideration of terms of the form $f(e_i\ot e_{i+sm}), f(e_i\ot
e_{i+sm-1})$ and $f(e_i\ot e_{i+sm+1})$ where $s \in {\mathbb Z}, s
\pgq 0$.

For any $s \pgq 0$ and $0 \ppq i \ppq m-1$, we have that $f(e_i\ot
e_{i+sm})$ is a linear combination of the elements $(a\oa)^k$,
$0\ppq k\ppq N-1$, and $(\oa a)^k$, $1\ppq k\ppq N$; $f(e_i\ot
e_{i+sm-1})$ is a linear combination of the $\oa(a\oa)^k$, $0\ppq
k\ppq N-1$; and $f(e_i\ot e_{i+sm+1})$ is a linear combination of
the $a(\oa a)^k$, $0\ppq k\ppq N-1$. For $m\pgq 3$, we give here a
list of the terms which must be considered in giving a possible
non-zero image of $f$; these are used without further comment
throughout the paper. (We assume $m \pgq 3$; the cases $m=1$ and
$m=2$ are different and are considered separately later.) For each
$i = 0, 1, \ldots , m-1$:
$$\begin{array}{llll}
{\bf m \mbox{\bf \ even, } m\neq 2}&t\mbox{ even }& f(e_i\ot e_{i+\alpha m }),&-p\ppq \alpha\ppq p, \\
&t\mbox{ odd, }t\neq m-1& f(e_i\ot e_{i+ \beta m-1}),&-p\ppq \beta\ppq p,\\
&& f(e_i\ot e_{i+ \gamma m+1}),&-p\ppq \gamma\ppq p,\\
&t=m-1& f(e_i\ot e_{i+ \beta m-1}),&-p\ppq \beta \ppq p+1,\\
&& f(e_i\ot e_{i+ \gamma m+1}),&-p-1\ppq \gamma\ppq p,\\
{\bf m \mbox{\bf \ odd, } m\neq 1}&t\mbox{ even, }t\neq m-1& f(e_i\ot e_{i+ (p-2\alpha)m}),&
  0\ppq \alpha\ppq p,\\
&& f(e_i\ot e_{i+ (p-2\beta-1)m-1 }),&0\ppq \beta\ppq p-1,\\
&& f(e_i\ot e_{i+ (p-2\gamma -1)m+1 }),&0\ppq\gamma\ppq p-1,\\
&t\mbox{ odd }& f(e_i\ot e_{i+ (p-2\alpha -1)m}),&0\ppq \alpha\ppq p-1,\\
&& f(e_i\ot e_{i+ (p-2\beta)m-1}),&0\ppq\beta\ppq p,\\
&& f(e_i\ot e_{i+ (p-2\gamma)m+1}),&0\ppq \gamma\ppq p,\\
&t=m-1& f(e_i\ot e_{i+ (p-2\alpha)m}),&0\ppq \alpha\ppq p,\\
&& f(e_i\ot e_{i+ (p-2\beta-1)m-1 }),&-1\ppq \beta\ppq p-1,\\
&& f(e_i\ot e_{i+ (p-2\gamma -1)m+1 }),&0\ppq \gamma\ppq p.
\end{array}$$

\subsection{Kernel of $\partial^{pm+t+1}:\Hom(P^{pm+t},\lan)\rightarrow \Hom(P^{pm+t+1},\lan)$.}

\sloppy An element $f$ in $\Hom(P^{pm+t},\lan)$ is determined by
$f(e_i\ot e_{i+n-2r})$ for appropriate $i\in\{0, 1, \ldots , m-1\}.$
Each term $f(e_i\ot e_{i+n-2r})$ can be written as a linear
combination of the basis elements in $e_i\lan e_{i+n-2r}$, and we
need to find the conditions on the coefficients in this linear
combination when $f$ is in $\ker\partial^{pm+t+1}$. To do this, we
separate the study into several cases: when $m$ is even, we consider
the three cases (i) $t=m-1$, (ii) $t$ even, and (iii) $t$ odd with
$t\neq m-1$, and, when $m$ is odd, we look at the four cases (i)
$t=m-1$, (ii) $t=m-2$, (iii) $t$ even with $t\neq m-1$, and (iv) $t$
odd with $t\neq m-2$. Note that in some cases a factor
 $N$ appears so that the results depend on the characteristic of
$\kk$. Note also that, when $m$ is odd, the case $\car\kk = 2$ needs
to be considered separately.

For $f\in\Hom(P^{pm+t},\lan)$, we write $\begin{cases}
f(e_i\ot e_{i+\alpha m})=\sum_{k=0}^{N-1}\sigma_{i,k}^{\alpha}(a_i\oa_i)^k+
\sum_{k=1}^N\tau_{i,k}^{\alpha}(\oa_{i-1}a_{i-1})^k,\\
f(e_i\ot e_{i+\beta m-1})=\sum_{k=0}^{N-1}\lambda_{i,k}^{\beta}(\oa a)^k\oa _{i-1},\\
f(e_i\ot e_{i+\gamma
m+1})=\sum_{k=0}^{N-1}\mu_{i,k}^{\gamma}(a\oa)^ka_i,
\end{cases}
$\\ with coefficients $\sigma_{i,k}^{\alpha},$
$\tau_{i,k}^{\alpha},$ $\lambda_{i,k}^\beta$ and $\mu_{i,k}^\gamma$
in $\kk$, for $i = 0, 1, \ldots , m-1$.

We shall now do two of the cases above, in order to illustrate the
general method. In both these cases, for consistency with our
convention on the writing of elements $e_i\ot e_{i+k}$ with
$i\in\set{0,1,\ldots,m-1}$ and $k\in\zz$, and to simplify, we shall
use the following notation: if $n\in\zz,$ $\rep{n}$ denotes the
representative of $n$ modulo $m$ in $\set{0,1,\ldots,m-1}.$ In
particular, $\rep{m}=\rep{0}$, and this must be borne in mind when
computing kernels and images. We write $\partial f$ (rather than
$\partial^{pm+t+1}f$) for the image of $f$ under the map
$\partial^{pm+t+1} \colon \Hom(P^{pm+t}, \L_N) \to \Hom(P^{pm+t+1},
\L_N)$.

\medskip

\noindent\textbf{Case $m$ even and $t$ even.} In this case, $f$ is
entirely determined by the $f(e_i\ot e_{i+\alpha'm})$ for
$-p\ppq{\alpha'}\ppq p$ and $i\in\set{0,1,\ldots,m-1}$ and $\partial
f$ is determined by\\
$\begin{cases}
\partial f(e_i\ot e_{i+\beta m-1}) \mbox{ for }  -p\ppq\beta \ppq p \mbox{ (or $-p\ppq \beta \ppq p+1$ if $t=m-2$)}\\
\partial f(e_i\ot e_{i+\gamma m+1})  \mbox{ for } -p\ppq\gamma\ppq p \mbox{ (or $-p-1\ppq \gamma \ppq p$ if $t=m-2$)}
\end{cases}$
with  $0\ppq i\ppq m-1.$  So $f$ is in the kernel of $\partial$ if
and only if $\partial f(e_i\ot e_{i+\beta m-1})=0=\partial f(e_i\ot
e_{i+\gamma m+1})$ for all $i,$ $\beta$ and $\gamma.$ Now
$\begin{array}{l}
\partial f(e_i\ot e_{i+\beta m-1})=\\
\ \ \ \begin{cases}
\left[(-1)^{(p-\beta)\frac{m}{2}}\sigma_{\rep{i-1},0}^\beta-\sigma_{i,0}^\beta\right]\oa(a\oa)^{N-1}
&\mbox{if $\beta>0$}\\
\left[(-1)^{(p-\beta)\frac{m}{2}}\sigma_{\rep{i-1},0}^\beta-\sigma_{i,0}^\beta\right]\oa+\sum_{k=1}^{N-1}\left[(-1)^{(p-\beta)\frac{m}{2}}\sigma_{\rep{i-1},k}^\beta-\tau_{i,k}^\beta\right](\oa
a)^k\oa&\mbox{if $\beta\ppq 0$}
\end{cases}\\
\partial f(e_i\ot e_{i+\gamma m+1})=\\
\ \ \ \begin{cases}
\left[\sigma_{i,0}^\gamma-(-1)^{(p-\gamma)\frac{m}{2}}\sigma_{\rep{i+1},0}^\gamma\right]a+\sum_{k=1}^{N-1}\left[\sigma_{i,k}^\gamma-(-1)^{(p-\gamma)\frac{m}{2}}\tau_{\rep{i+1},k}^\gamma\right]a(\oa
a)^k&\mbox{if $\gamma\pgq 0$}\\
\left[\sigma_{i,0}^\gamma-(-1)^{(p-\gamma)\frac{m}{2}}\sigma_{\rep{i+1},0}^\gamma\right]a(\oa
a)^{N-1}&\mbox{if $\gamma>0$.}
\end{cases}
\end{array}
$

For any $\alpha'$ with $-p\ppq \alpha'\ppq p$, we can see that
$\sigma_{0,0}^{\alpha'}$ determines the other
$\sigma_{i,0}^{\alpha'}$ (setting $\beta$ or $\gamma$ to be
$\alpha'$). Combining the cases $\beta\ppq 0$ and $\gamma\pgq 0$
tells us that the $\tau_{i,k}^{\alpha'}$ for $k=1, \ldots , N-1$ are
determined by the $\sigma_{i,k}^{\alpha'}.$ Finally the
$\tau_{i,N}^{\alpha'}$ are arbitrary, and the
$\sigma_{i,k}^{\alpha'}$ for $k=1, \ldots , N-1$ are arbitrary.
Therefore the dimension of $\ker\partial^{pm+t+1}$ in this case is
$(2p+1)+ (2p+1)m + (2p+1)m(N-1).$

\medskip

\noindent\textbf{Case $m$ even and $t$ odd, $t\neq m-1$.} In this
case, $f$ is entirely determined by $\begin{cases}
f(e_i\ot e_{i+\beta'm-1}) \mbox{ for } -p\ppq \beta'\ppq p,\\
f(e_i\ot e_{i+\gamma'm+1})\mbox{ for } -p\ppq \gamma'\ppq p,
\end{cases}$ with $i\in \set{0,1,\ldots,m-1}$ \\
and $\partial f\in \Hom(P^{(p+1)m},\lan)$ is determined by $\partial
f(e_i\ot e_{i+\alpha m})$ for $-p\ppq \alpha\ppq p.$ Now $$\partial
f(e_i\ot e_{i+\alpha m})=
\begin{cases}
N\left[\lambda_{i,0}^0+(-1)^{(p+1)\frac{m}{2}}\lambda_{\rep{i+1},0}^0+(-1)^{(p+1)\frac{m}{2}}\mu_{\rep{i-1},0}^0+\mu_{i,0}^0\right](a\oa)^N\mbox{
  if
$\alpha=0$}\\
\left[\lambda_{i,N-1}^\alpha+(-1)^{(p+1+\alpha)\frac{m}{2}}\lambda_{\rep{i+1},N-1}^\alpha+(-1)^{(p+1+\alpha)\frac{m}{2}}\mu_{\rep{i-1},0}^\alpha+\mu_{i,0}^\alpha\right](a\oa)^N
\\\hspace*{.5cm}+ \sum_{k=0}^{N-2}\lambda_{i,k}^\alpha (\oa
  a)^{k+1}+(-1)^{(p+1+\alpha)\frac{m}{2}}\sum_{k=0}^{N-2}\lambda_{\rep{i+1},k}^\alpha (a\oa)^{k+1}\mbox{ if $\alpha>0$}\\
\left[\lambda_{i,0}^\alpha+(-1)^{(p+1+\alpha)\frac{m}{2}}\lambda_{\rep{i+1},0}^\alpha+(-1)^{(p+1+\alpha)\frac{m}{2}}\mu_{\rep{i-1},N-1}^\alpha+\mu_{i,N-1}^\alpha\right](a\oa)^N
\\\hspace*{.5cm}+
(-1)^{(p+1+\alpha)\frac{m}{2}}\sum_{k=0}^{N-2}\mu_{\rep{i-1},k}^\alpha
  (\oa a)^{k+1}+\sum_{k=0}^{N-2}\mu_{i,k}^\alpha
  (a\oa)^{k+1}\mbox{ if $\alpha<0$.}
\end{cases}
$$

So for $\alpha=0,$ $\lambda_{i,0}^0, $ for $ i \in
\set{0,\ldots,m-1}$ and $\mu_{0,0}^0$ determine the other
$\mu_{i,0}^0,$ unless the characteristic of $\kk$ divides $N.$ For
$\alpha>0,$ $\mu_{i,0}^\alpha$ and $\lambda_{0,N-1}^\alpha$
determine the other $\lambda_{i,N-1}^\alpha,$ and for $\alpha<0,$
$\lambda_{i,0}^\alpha$ and $\mu_{0,N-1}^\alpha$ determine the other
$\mu_{i,N-1}^\alpha.$ Moreover, for $\alpha>0$ the
$\lambda_{i,k}^\alpha$ for $k=0, \ldots , N-2$ are all 0, and for
$\alpha<0$ the $\mu_{i,k}^\alpha$ for $k=0, \ldots , N-2$ are 0. The
remaining $\lambda_{i,k}^\alpha,$ $\mu_{i,k}^\alpha$ are arbitrary,
and there are $2m(N-1)+pm(N-1)+pm(N-1) = 2(p+1)m(N-1)$ of these
coefficients. This gives dimension: $(m+1)+p(m+1)+p(m+1)+
2(p+1)m(N-1)$ if $\car\kk\nmid N$ and
$(m+1)+p(m+1)+p(m+1)+2(p+1)m(N-1)+m-1$ if $\car\kk\mid N.$

\medskip

Note that when computing the dimensions in the case $m$ odd and
$\car\kk\neq 2$, we often need to introduce additional cases
depending on the value of $p$ modulo $4$ and the value of $m$ modulo
4. We summarise the results in the following two propositions, which
separate the cases $\car\kk\nmid N$ and $\car\kk\mid N$.

\bprop\label{dimensionsker} If $\car\kk\nmid N,$ the dimension of
$\ker\partial^{pm+t+1}$ is as follows.
$$\begin{array}{lll}
{\mbox{\textbf{  If $m$ even }}} & (2p+1)(mN+1)& \mbox{ if $t$ even,}\\
 & (2p+1)(mN+1)+m(N-1)& \mbox{ if $t$ odd, $t\neq m-1$ }\\
 & (2p+1)(mN+1)+m(N-1)+2 & \mbox{ if $t=m-1.$}\\\\
\mbox{\textbf{  If $m$ odd  }} & (2p+1)(mN+1)& \mbox{ if $t\neq m-1,$ $p+t$ even,}\\
\mbox{\textbf{  and $\car\kk=2$ }} & (2p+1)(mN+1)+m(N-1)& \mbox{ if $t\neq m-1,$ $p+t$ odd,}\\
 & (2p+1)(mN+1)+ 2& \mbox{ if $t=m-1,$ $p$ even,}\\
 & (2p+1)(mN+1)+m(N-1)+ 2& \mbox{ if $t=m-1,$ $p$ odd.}\\\\
\mbox{\textbf{  If $m$ odd  }} & (2p+1)mN+\frac{p}{2}+1& \mbox{ if $t\md{0}$, $t\neq m-1, $ $p$ even,}\\
\mbox{\textbf{  and $\car\kk\neq 2$ }} & (2p+1)mN+\frac{p}{2}& \mbox{ if $t\md{2}$, $t\neq m-1, $ $p$ even,}\\
 & (2p+1)mN+m(N-1)+\frac{p+1}{2}& \mbox{ if $m+t\md{1}$, $t\neq m-1, $ $p$ odd, }\\
 & (2p+1)mN+m(N-1)+\frac{p-1}{2}& \mbox{ if $m+t\md{3}$, $t\neq m-1, $ $p$ odd, }\\
 & (2p+1)mN+m(N-1)+  \frac{p}{2}+1& \mbox{ if $t\md{1}$, $p$ even, }\\
 & (2p+1)mN+m(N-1)+ \frac{p}{2}& \mbox{ if $t\md{3}$, $p$ even, }\\
 & (2p+1)mN+\frac{p-1}{2}& \mbox{ if $t\md{m}$, $p$ odd, }\\
 & (2p+1)mN+\frac{p+1}{2}& \mbox{ if $t\md{m+2}$, $p$ odd, }\\
 & (2p+1)mN+\frac{p}{2}+3& \mbox{ if $t=m-1\md{0},$ $p$ even, }\\
 & (2p+1)mN+\frac{p}{2}+2& \mbox{ if $t=m-1\md{2},$ $p$ even, }\\
 & (2p+1)mN+m(N-1)+\frac{p+1}{2}& \mbox{ if $t=m-1,$ $p$ odd.}
\end{array}
$$\eprop

\bprop\label{dimensionsker2} If $\car\kk\mid N,$ the dimension of
$\ker\partial^{pm+t+1}$ is given by the following.
$$
\begin{array}{lll}
{\mbox{\textbf{ If $m$ even }}}
& (2p+1)(mN+1) & \mbox{ if $t$ even,}\\
& (2p+2)mN+2p  & \mbox{ if $t$ odd, $t\neq m-1$ }\\
& (2p+2)(mN+1) & \mbox{ if $t=m-1.$}\\\\
\mbox{\textbf{  If $m$ odd } }
&  (2p+1)(mN+1)   & \mbox{ if $t\neq m-1,$ $p+t$ even,}\\
\mbox{\textbf{  and $\car\kk = 2$ }}
&  (2p+2)mN+2p    & \mbox{ if $t\neq m-1,$ $p+t$ odd,}\\
&  (2p+1)(mN+1)+2 & \mbox{ if $t=m-1,$ $p$ even,}\\
&  (2p+2)(mN+1)   & \mbox{ if $t=m-1,$ $p$ odd.}\\\\
\mbox{\textbf{  If $m$ odd  }}
& (2p+1)mN+ \frac{p}{2}+1 & \mbox{ if $t\md{0}$, $p$ even, }\\
\mbox{\textbf{  and $\car\kk\neq 2$} }
& (2p+1)mN+ \frac{p}{2}   & \mbox{ if $t\md{2}$, $p$ even, }\\
& (2p+2)mN+ \frac{p-1}{2} & \mbox{ if $t$ even, $p\md{1}$}\\
& (2p+2)mN+ \frac{p+1}{2} & \mbox{ if $m+t\md{1}$, $p\md{3}$}\\
& (2p+2)mN+ \frac{p-3}{2} & \mbox{ if $m+t\md{3}$, $p\md{3}$}\\
& (2p+2)mN+ \frac{p}{2}   & \mbox{ if $t$ odd, $p \md{0}$ }\\
& (2p+2)mN+ \frac{p}{2}+1 & \mbox{ if $t\md{1}$, $p\md{2}$}\\
& (2p+2)mN+ \frac{p}{2}-1 & \mbox{ if $t\md{3}$, $p\md{2}$}\\
& (2p+1)mN+ \frac{p-1}{2} & \mbox{ if $t\md{m}$, $p$ odd}\\
& (2p+1)mN+ \frac{p+1}{2} & \mbox{ if $t\md{m+2}$, $p$ odd}\\
& (2p+2)mN+ \frac{p+1}{2} & \mbox{ if $t=m-1$, $p\md{3}$}\\
& (2p+2)mN+ \frac{p-1}{2} & \mbox{ if $t=m-1$, $p\md{1}$}\\
& (2p+1)mN+ \frac{p}{2}+3 & \mbox{ if $t=m-1\md{0}$, $p$ even}\\
& (2p+1)mN+ \frac{p}{2}+2 & \mbox{ if $t=m-1\md{2}$, $p$ even.}
\end{array}
$$\eprop

\subsection{Image of $\partial^{pm+t}:\Hom(P^{pm+t-1},\lan)\rightarrow \Hom(P^{pm+t},\lan)$.}

We explain here how to compute $\im\partial^{pm+t}$ explicitly
(rather than simply compute its dimension with the rank and nullity
theorem) on the same cases as for the kernels, in order to be able
to give a basis of cocycle representatives for each Hochschild
cohomology group $\HH^n(\lan)$. For an element
$g\in\Hom(P^{pm+t-1},\lan),$ we consider $\partial^{pm+t}g$, and
find the dimension of the subspace of $\ker(\partial^{pm+t+1})$
spanned by all the $\partial^{pm+t}g$. We shall again use the
notation $\rep{n}$ for the representative modulo $m$ of the integer
$n$, and write $\partial g$ for $\partial^{pm+t}g$.

\medskip
\noindent\textbf{Case $m$ even and $t$ even, $t\neq 0$.}

The map $\partial g$ is determined by $\partial g(e_i\ot e_{i+\alpha
m})$ for $-p\ppq \alpha \ppq p$ and $g$ is determined by
$\begin{cases}
g(e_i\ot e_{i+ \beta' m-1}) \mbox{ for $-p\ppq \beta' \ppq p$}\\
g(e_i\ot e_{i+ \gamma' m+1}) \mbox{ for $-p\ppq \gamma '\ppq p$}
\end{cases}$ with $i\in\set{0,1,\ldots,m-1}.$ Now $$
\partial g(e_i\ot e_{i+\alpha m})=\begin{cases}
N\left[(-1)^{p\frac{m}{2}}\mu_{\rep{i-1},0}^0+\mu_{i,0}^0+\lambda_{i,0}^0+(-1)^{p\frac{m}{2}}\lambda_{\rep{i+1},0}^0\right](a\oa)^N      \mbox{ if $\alpha =0$}\\
\left[(-1)^{(p+\alpha)\frac{m}{2}}\mu_{\rep{i-1},0}^\alpha+\mu_{i,0}^\alpha+\lambda_{i,N-1}^\alpha+(-1)^{(p+\alpha)\frac{m}{2}}\lambda_{\rep{i+1},N-1}^\alpha\right](a\oa)^N   \\\hspace{1cm}+\sum_{k=0}^{N-2}\lambda_{i,k}^\alpha(\oa a)^{k+1}+\sum_{k=0}^{N-2}\lambda_{\rep{i+1},k}^\alpha(a\oa )^{k+1}     \mbox{ if $\alpha >0$}\\
\left[(-1)^{(p+\alpha)\frac{m}{2}}\mu_{\rep{i-1},N-1}^\alpha+\mu_{i,N-1}^\alpha+\lambda_{i,0}^\alpha+(-1)^{(p+\alpha)\frac{m}{2}}\lambda_{\rep{i+1},0}^\alpha\right](a\oa)^N
\\
\hspace{1cm}+\sum_{k=0}^{N-2}\mu_{i,k}^\alpha(a\oa
)^{k+1}+(-1)^{(p+\alpha)\frac{m}{2}}\sum_{k=0}^{N-2}\mu_{\rep{i-1},k}^\alpha(\oa
a)^{k+1}        \mbox{ if $\alpha <0$.}
\end{cases}
$$

We can check that in each case, the sum or the alternate sum over
$i\in\set{0,\ldots,m-1}$ of the coefficients of the $(a\oa)^N$ is
zero (recall that $\rep{m}=\rep{0}$ and that $m$ is even), so for
each $\alpha$ we add $m-1$ to the dimension, except when
$\car\kk\mid N$ and $\alpha=0$, in which case, $\partial
g(e_i\otimes e_i) = 0$ and this makes no contribution to
$\dim\im\partial^{pm+t}$. Moreover, for $\alpha\neq 0$, we can take
arbitrary coefficients for the $(a\oa)^{k+1}$ with $0\ppq k\ppq
N-2$, and the coefficients of the $(\oa a)^{k+1}$ are completely
determined by those of the $(a\oa)^{k+1}$.

So $\dim_\kk\im\partial^{pm+t}=
\begin{cases}
(2p+1)(m-1)+2pm(N-1)&\mbox{if $\car\kk\nmid N$}\\
2p(mN-1)&\mbox{if $\car\kk\mid N$}.
\end{cases}$

\medskip

\noindent\textbf{Case $m$ even and $t$ odd, $t\neq m-1$.}

The element $\partial g$ is determined by $\begin{cases}\partial
g(e_i\ot e_{i+\beta m-1})\mbox{ for $-p\ppq \beta \ppq p$ }\\
\partial g(e_i\ot e_{i+\gamma m+1}) \mbox{ for $-p\ppq\gamma\ppq p$}\end{cases}$
and $g$ is determined by $g(e_i\ot e_{i+\alpha'm})$ for $-p\ppq
\alpha' \ppq p.$ Now
$$\begin{array}{l}
\partial g(e_i\ot e_{i+\beta m-1})=\\\
\begin{cases}
 -\left[(-1)^{(p+\beta+1)\frac{m}{2}}\sigma_{\rep{i-1},0}^\beta+\sigma_{i,0}^\beta\right](\oa a)^{N-1}\oa
 &\mbox{ if $\beta>0$}\\
-\left[(-1)^{(p+\beta+1)\frac{m}{2}}\sigma_{\rep{i-1},0}^\beta+\sigma_{i,0}^\beta\right]\oa-\sum_{k=1}^{N-1}\left[(-1)^{(p+\beta+1)\frac{m}{2}}\sigma_{\rep{i-1},k}^\beta+\tau_{i,k}^\beta\right](\oa
a)^k\oa &\mbox{ if $\beta\ppq 0$}
\end{cases}\\
\partial g(e_i\ot e_{i+\gamma m+1})=\\\
\begin{cases}
\left[\sigma_{i,0}^\gamma+(-1)^{(p+\gamma+1)\frac{m}{2}}\sigma_{\rep{i+1},0}^\gamma\right]a+
\sum_{k=1}^{N-1}\left[\sigma_{i,k}^\gamma+(-1)^{(p+\gamma+1)\frac{m}{2}}\tau_{\rep{i+1},k}^\gamma\right]a(\oa a)^k &\mbox{ if $\gamma \pgq 0$}\\
\left[\sigma_{i,0}^\gamma+(-1)^{(p+\gamma+1)\frac{m}{2}}\sigma_{\rep{i+1},0}^\gamma\right]a(\oa
a)^{N-1} &\mbox{ if $\gamma <0$.}
\end{cases}
\end{array}$$

For $0 < \alpha' \ppq p$ (resp. $-p \ppq \alpha' < 0$), the
coefficient of $(\oa a)^{N-1}\oa$ (resp. $\oa$) in $\partial g(e_i
\otimes e_{i+\alpha'm-1})$ is equal, up to sign, to the coefficient
of $a$ (resp. $a(\oa a)^{N-1}$) in $\partial g(e_i\otimes
e_{i+\alpha' m+1})$. The coefficient of $\oa$ in $\partial
g(e_i\otimes e_{i-1})$ (that is, when $\alpha'=0$) is equal, up to
sign, to that of $a$ in $\partial g(e_i\otimes e_{i+1})$. Moreover,
either the sum or the alternate sum over $i\in\set{0,\ldots,m-1}$ of
these coefficients is zero (depending on the parity of
$(p+\alpha'+1)\frac{m}{2}$). They contribute $m-1$ to the dimension
of $\im\partial^{pm+t}$ for each $\alpha'$. The remaining
coefficients of the $(\oa a)^k\oa$ and $a(\oa a)^k$, for $k = 1,
\ldots , N-1$, are arbitrary except for the case $\gamma = 0$, where
the coefficient of $a(\oa a)^k$ is equal, up to sign, to that of
$(\oa a)^k\oa$ in $\partial g(e_i\otimes e_{i-1})$ (for $\beta =
0$).

So $\dim_\kk\im\partial^{pm+t}=(2p+1)(m-1)+(2p+1)m(N-1).$

\subsection{Dimension of $\HH^n(\lan)$.}

We will now give the dimension of each of the Hochschild cohomology
groups of $\L_N$. Note that the dimensions can be verified using
Propositions \ref{dimensionsker} and \ref{dimensionsker2}, together
with the rank and nullity theorem in the form
$\dim\im(\partial^{pm+t})=\dim\Hom(P^{pm+t-1},\lan)-\dim\ker(\partial^{pm+t})$.

\bprop\label{prop:dim} If $m$ is even, or if $m$ is odd and
$\car\kk=2$, then $\dim_\kk\HH^0(\lan)=Nm+1$ and, for $pm+t \neq 0$,
$$\dim_\kk\HH^{pm+t}(\lan)=\begin{cases}
4p+2+m(N-1) &\mbox{ if $t\neq m-1,$ $\car\kk\nmid N$,}\\
4p+4+m(N-1) &\mbox{ if $t=m-1$,  $\car\kk\nmid N$,}\\
4p+1+mN&\mbox{ if $t\neq m-1,$ $\car\kk\mid N$,}\\
4p+3+mN&\mbox{ if $t=m-1$,  $\car\kk\mid N$.}
\end{cases}
$$ \eprop

\bprop  If $m$ is odd, $\car\kk\neq 2$ and $\car\kk\nmid N$, then
$\dim_\kk\HH^0(\lan)=Nm+1$ and, for $pm+t \neq 0$,
$$\dim_\kk\HH^{pm+t}(\lan)=m(N-1)+p+
\begin{cases}
1 &\mbox{ if $t=0$, $p$ even,}\\
3 &\mbox{ if $t=0$, $p$ odd, $\frac{m-1}{2}$ even,}\\
1 &\mbox{ if $t=0$, $p$ odd, $\frac{m-1}{2}$ odd, }\\
1 &\mbox{ if $t$ even,  $t\neq 0,m-1,$ $p$ even,}\\
-1 &\mbox{ if $m+t\md{3}$ even,  $t\neq 0,m-1,$ $p$ odd,}\\
1 &\mbox{ if $m+t\md{1}$ even,  $t\neq 0,m-1,$ $p$ odd, }\\
2   &\mbox{ if $t\md{1}$,  $p$ even, }\\
0 &\mbox{ if $t\md{3}$,  $p$ even, }\\
0   &\mbox{ if $t$ odd,  $p$ odd, }\\
3 &\mbox{ if $t=m-1$, $p$ even,}\\
1 &\mbox{ if $t=m-1$, $p$ odd.}
\end{cases}
$$\eprop

\bprop  If $m$ is odd, $\car\kk\neq 2$ and $\car\kk\mid N$, then
$\dim_\kk\HH^0(\lan)=Nm+1$ and, for $pm+t \neq 0$,
$$\dim_\kk\HH^{pm+t}(\lan)=mN+p+
\begin{cases}
1 & \mbox{if $t=0$, $p\md{0,1}$ }\\
0 & \mbox{if $t=0$, $p\md{2}$ }\\
3 & \mbox{if $t=0$, $p\md{3}$, $m\md{1}$}\\
0 & \mbox{if $t=0$, $p\md{3}$, $m\md{3}$}\\
1 & \mbox{if $t$ even, $t\neq 0,m-1$, $p+t\md{0}$}\\
0   & \mbox{if $t$ even, $t\neq 0,m-1$, $p+t\md{2}$}\\
0   & \mbox{if $m+t\md{1}$, $t\neq 0,m-1$, $p\md{1}$, }\\
-1 & \mbox{if $m+t\md{3}$, $t\neq 0,m-1$, $p\md{1}$}\\
1 & \mbox{if $m+t\md{1}$, $t\neq 0,m-1$, $p\md{3}$}\\
-2 & \mbox{if $m+t\md{3}$, $t\neq 0,m-1$, $p\md{3}$}\\
1 & \mbox{if $t\md{1}$, $p\md{0}$}\\
0   & \mbox{if $t\md{3}$, $p\md{0}$}\\
2 & \mbox{if $t\md{1}$, $p\md{2}$}\\
-1 & \mbox{if $t\md{3}$, $p\md{2}$}\\
-1 & \mbox{if $t$ odd, $p+m-t\md{1}$}\\
0   & \mbox{if $t$ odd, $p+m-t\md{3}$}\\
3 & \mbox{if $t=m-1$, $p+m\md{1}$}\\
2 & \mbox{if $t=m-1$, $p+m\md{3}$}\\
0   & \mbox{if $t=m-1$, $p\md{1}$}\\
1 & \mbox{if $t=m-1$, $p\md{3}$.}
\end{cases}
$$\eprop

\bigskip

\section{The centre of the algebra $\L_N$}

We start by describing $Z(\L_N)$, the centre of $\L_N$, since it is
well known that $\HH^0(\L_N) = Z(\L_N)$. Note that the next result
requires $m\pgq 2$; the case $m=1$ is dealt with separately in
Theorem \ref{thm:m1centre}.

\bt\label{thm:centre} Suppose that $N \pgq 1, m \pgq 2$. Then
$\dim\HH^0(\L_N) = Nm+1$ and the set
$$\set{1, (a_i\oa_i)^N , [(a_i\oa_i)^s+(\oa_ia_i)^s] \mbox{ for } i=0, 1, \ldots, m-1 \mbox{
and } s = 1, \ldots , N-1}$$ is a $\kk$-basis of $\HH^0(\lan)$.

Moreover, $\HH^0(\lan)$ is generated \textit{as an algebra} by the
set $$\begin{array}{ll}
\set{1, a_i\oa_i \mbox{ for } i = 0, 1 \ldots , m-1} & \mbox{if $N=1$;}\\
\set{1, (a_i\oa_i)^N, [a_i\oa_i+\oa_ia_i] \mbox{ for } i = 0, 1
\ldots , m-1} & \mbox{if $N>1$.}
\end{array}$$ \et

The proof is straightforward since, for $N>1$, we have
$[(a_i\oa_i)^s+(\oa_ia_i)^s] = ([a_i\oa_i+\oa_ia_i])^s$.

\bigskip

\section{The Hochschild cohomology ring $\HH^*(\L_N)$ for $m \pgq 3$ and $m$ even.}

In this section, we assume that $N \pgq 1$, $m \pgq 3$ and $m$ even
(so that, necessarily, $m\pgq 4$) and give all the details of
$\HH^*(\L_N)$. In the subsequent sections we will consider the case
$N \pgq 1$, $m \pgq 3$ and $m$ odd, and leave many of the details to
the reader since the computations are similar.

\subsection{Basis of $\HH^n(\L_N)$}

For each $n \pgq 1$, we describe the elements of a basis of the
Hochschild cohomology group $\HH^n(\L_N)$ in terms of cocycles in
$\Hom(P^{n},\lan)$, that is, we give a set of elements in
$\ker\partial^{n+1}$, each of which is written as a map in
$\Hom(P^{n},\lan)$, such that the corresponding set of cosets in
$\ker\partial^{n+1}/\im\partial^{n}$ with these representatives
forms a basis of $\ker\partial^{n+1}/\im\partial^{n} = \HH^n(\L_N)$.
Following standard usage, our notation does not distinguish between
an element of $\HH^n(\L_N)$ and a cocycle in $\Hom(P^{n},\lan)$
which represents that element of $\HH^n(\L_N)$. However the precise
meaning is always clear from the context.  For each cocycle in
$\Hom(P^{n},\lan)$, we simply write the images of the generators
$e_i\ot e_{i+n-2r}$ in $P^{n}$ which are non-zero. We keep to these
conventions throughout the whole paper.

\bprop\label{prop:basisHHn} Suppose that $N \pgq 1$, $m \pgq 3$ and
$m$ is even. For each $n \pgq 1$, the following elements define a
basis of $\HH^n(\L_N)$.
\begin{enumerate}
\item For $n$ even, $n \pgq 2$:
\begin{enumerate}
\item For $-p\ppq \alpha \ppq p$: $\chi_{n,\alpha}:\ e_i\ot e_{i+\alpha m}\mapsto
(-1)^{\left(\frac{n}{2}-\alpha\frac{m}{2}\right)i}e_i$ for $i = 0,
1, \ldots , m-1$;
\item For $-p\ppq \alpha \ppq p$: $\pi_{n,\alpha}:\ e_0\ot e_{\alpha m}\mapsto
(a_0\oa_0)^N$;
\item For each $j= 0,
1, \ldots , m-2$ and each $s = 1, \ldots , N-1$:\\ $F_{n,j,s}:
\begin{cases}
e_j\ot e_j\mapsto (a_j\oa_j)^s\\
e_{j+1}\ot e_{j+1}\mapsto (-1)^{\frac{n}{2}}(\oa_ja_j)^s;
\end{cases}$
\item For each $s = 1, \ldots , N-1$:\\ $F_{n,m-1,s}:
\begin{cases}
e_{m-1}\ot e_{m-1}\mapsto (a_{m-1}\oa_{m-1})^s\\
e_0\ot e_0\mapsto (-1)^{\frac{n}{2}}(\oa_0a_0)^s;
\end{cases}$
\item Additionally in the case $\car\kk\mid N$, for each $j= 1, \ldots , m-1$: $\theta_{n,j}:e_j\ot
e_j\mapsto (a_j\oa_j)^N$.
\end{enumerate}
\item For $n$ odd, $n \pgq 1$:
\begin{enumerate}
\item For $-p\ppq \gamma<0$ and, if $t=m-1$, $\gamma=-p-1$ also:\\
\hspace*{.5cm} $\varphi_{n,\gamma}:\ e_i\ot e_{i+\gamma m+1}\mapsto
(-1)^{\left(\frac{n-1-\gamma m}{2}\right)i}a_{i}(\oa_i a_i)^{N-1}$
for all $i= 0, 1, \ldots , m-1$;
\item For $0\ppq \gamma \ppq p$:
$\varphi_{n,\gamma}:\ e_i\ot e_{i+\gamma m+1}\mapsto
(-1)^{\left(\frac{n-1-\gamma m}{2}\right)i}a_{i}$ for all $i= 0, 1,
\ldots , m-1$;
\item For $0<\beta\ppq p$ and, if $t=m-1,$ $\beta=p+1$ also:\\
\hspace*{.5cm} $\psi_{n,\beta}:e_i\ot e_{i+\beta m-1}\mapsto
(-1)^{\frac{n-1-\beta m}{2}i}\oa_{i-1}(a_{i-1}\oa_{i-1})^{N-1}$ for
all $i= 0, 1, \ldots , m-1$;
\item For $-p\ppq \beta\ppq 0$:
$\psi_{n,\beta}:e_i\ot e_{i+\beta m-1}\mapsto (-1)^{\frac{n-1-\beta
m}{2}i}\oa_{i-1}$ for all $i = 0, 1, \ldots , m-1$;
\item For each $j= 0,
1, \ldots , m-1$ and each $s=1, \ldots , N-1$: $E_{n,j,s}:e_j\ot
e_{j+1}\mapsto (a_j\oa_j)^sa_j$;
\item Additionally in the case $\car\kk\mid N$, for each $j= 1, \ldots , m-1$:
$\omega_{n,j}:e_j\ot e_{j+1}\mapsto a_j$.
\end{enumerate}
\end{enumerate}
\eprop

\bigskip

\subsection{Computing products in cohomology}

In order to compute the products in the Hochschild cohomology ring,
we require a lifting of each element in the basis of $\HH^n(\L_N)$,
that is, for each cocycle $f\in\Hom(P^{n},\lan)$, we give a map
${\mathcal L}^*f:P^{*+n}\rightarrow P^*$ such that
\begin{enumerate}[(1)]
\item the composition
$P^{n}\stackrel{{\mathcal
L}^0{f}}{\longrightarrow}P^0\stackrel{\partial^0}{\longrightarrow}\lan$
is equal to $f$, and
\item the square $\xymatrix{P^{q+n}\ar^{{\mathcal L}^q{f}}[r]\ar_{\partial^{q+n}}[d]
& P^q\ar^{\partial^q}[d]\\P^{q-1+n} \ar_{{\mathcal
L}^{q-1}{f}}[r]&P^{q-1}}$ commutes for every $q\pgq 1.$
\end{enumerate}
It should be noted that such liftings always exist but are not
unique. It is easy to check that the maps given in this paper are
indeed liftings of the appropriate cocycle.

\bigskip

For homogeneous elements $\eta \in \HH^n(\L_N)$ and $\theta \in
\HH^k(\L_N)$ represented by cocycles $\eta : P^n \to \L_N$ and
$\theta : P^k \to \L_N$ respectively, the cup product (or Yoneda
product) $\eta\theta \in \HH^{n+k}(\L_N)$ is the coset represented
by the composition $P^{n+k} \stackrel{{\mathcal
L}^n\theta}{\longrightarrow} P^n \stackrel{\eta}{\longrightarrow}
\L_N$. The fact that the cup product is well-defined means that the
composition $\eta\theta$ does not depend on the choice of
representative cocycles for $\eta$ and $\theta$ or on the choice of
the liftings of these cocycles.

\bigskip

\subsection{Liftings for $N = 1$}

We continue to assume that $m$ is even in the details that follow.

\bprop For $n \pgq 1$ and the algebra $\L_1$, $m \pgq3$ and $m$
even, the following maps are liftings of the cocycle basis of
$\HH^n(\L_1)$ given in Proposition \ref{prop:basisHHn}.
\begin{enumerate}
\item For $n$ even, $n\pgq 2$, and for $-p\ppq \alpha \ppq p$:
\begin{enumerate}[$\bullet$]
  \item ${\mathcal L}^q{\chi}_{n,\alpha}:\ e_i\ot e_{i+\alpha m+q-2\ell}\mapsto
  (-1)^{\left(\frac{n}{2}-\alpha\frac{m}{2}\right)i}e_i\ot e_{i +q-2\ell}$
  for all $i=0, 1, \ldots , m-1$ and all $0\ppq \ell\ppq q$,
  \item ${\mathcal L}^q{\pi}_{n,\alpha}:\ e_0\ot e_{\alpha m +q-2\ell}\mapsto
  a_0\oa_0e_0\ot e_{q-2\ell}$ for all $0\ppq \ell\ppq q$,
\end{enumerate}
\item For $n$ odd:
\begin{enumerate}[$\bullet$]
  \item For $\begin{cases}
            -p\ppq\gamma\ppq p & \mbox{ if $t\neq m-1$}\\
            -p-1\ppq \gamma\ppq p & \mbox{ if $t=m-1$:}
            \end{cases}$\\
  ${\mathcal L}^q{\varphi}_{n,\gamma}:\ e_i\ot e_{i+\gamma m+1 +q-2\ell}\mapsto
  (-1)^{\left(\frac{n-1}{2}-\gamma\frac{m}{2}\right)i+q}e_i\ot e_{i
  +q-2\ell}\;a_{i+q-2\ell}$\\
  for all $i=0, 1, \ldots , m-1$ and all $0\ppq \ell\ppq q$,
  \item For $\begin{cases}
            -p\ppq\beta\ppq p&\mbox{ if $t\neq m-1$}\\
            -p\ppq\beta\ppq p+1&\mbox{ if $t=m-1$:}
            \end{cases}$\\
  ${\mathcal L}^q{\psi}_{n,\beta}:\ e_i\ot e_{i+ \beta m-1+q-2\ell}\mapsto
  (-1)^{\left(\frac{n-1}{2}-\beta\frac{m}{2}\right)i} e_i\ot e_{i
  +q-2\ell}\;\oa_{i+q-2\ell-1}$\\
  for all $i=0, 1, \ldots , m-1$ and all $0\ppq \ell\ppq q$.
\end{enumerate}
\end{enumerate}
\eprop

\bigskip

\subsection{Liftings for $N > 1$}

\bprop\label{prop:liftings} For $n \pgq 1$ and the algebra $\L_N$
with $N>1, m \pgq3$ and $m$ even, the following maps give a lifting
of the cocycles in Proposition \ref{prop:basisHHn} above. Once
again, to simplify cases and for consistency with our conventions,
we use the notation $\rep{n}\in\set{0,1,\ldots,m-1}$ for the
representative of the integer $n$ modulo $m.$
\begin{enumerate}
\item For $n$ even, $n \pgq 2$:
\begin{enumerate}[$\bullet$]
\item For $-p\ppq \alpha \ppq p$:\\
${\mathcal L}^q{\chi}_{n,\alpha}(e_i\ot e_{i+q-2\ell+\alpha m})=\\
  \begin{cases}
  (-1)^{\frac{n-\alpha m}{2}i}e_i\ot e_{i+q-2\ell} &\mbox{if $(q-2\ell)\alpha\pgq 0$}\\
  0&\mbox{if $(q-2\ell)\alpha<0$ and $\abs{q-2\ell}>2$}\\
  (-1)^{\frac{n-\alpha m}{2}i} (a\oa)^{N-1}e_i\ot e_{i+2}(\oa a)^{N-1} &
  \mbox{if $\alpha<0$ and $q-2\ell=2$}\\
  (-1)^{\frac{n-\alpha m}{2}i} (\oa a)^{N-1}e_i\ot e_{i-2}(a\oa )^{N-1} &
  \mbox{if  $\alpha>0$ and $q-2\ell=-2$}\\
  (-1)^{\frac{n-\alpha m}{2}i} \left[\sum_{k=0}^{N-1}(a\oa)^ke_i\ot
  e_{i+1}(\oa a)^{N-k-1}\right.\\
  \hspace*{2cm}\left.+(-1)^{\frac{q+1}{2}}\sum_{k=0}^{N-2}(a\oa)^ka e_{\rep{i+1}}\ot e_{\rep{i+1}-1} a(\oa a)^{N-k-2}\right] &
  \mbox{if $\alpha<0$ and $q-2\ell=1$}\\
  (-1)^{\frac{n-\alpha m}{2}i} \left[\sum_{k=0}^{N-1}(\oa a)^ke_i\ot e_{i-1}(a\oa)^{N-k-1}\right.\\
  \hspace*{2cm}\left.+(-1)^{\frac{q+1}{2}}\sum_{k=0}^{N-2}(\oa a)^k\oa e_{\rep{i-1}}\ot e_{\rep{i-1}+1} \oa(a\oa)^{N-k-2}\right] &
  \mbox{if $\alpha>0$ and $q-2\ell=-1$}
  \end{cases}$ for all $i=0, 1, \ldots , m-1$ and all $0\ppq\ell\ppq q$.
\item For $-p\ppq \alpha \ppq p$:\\
${\mathcal L}^q{\pi}_{n,\alpha}(e_0\ot e_{q-2\ell+\alpha m})=
  \begin{cases}
  (a\oa)^N e_0\ot e_{q-2\ell} & \mbox{if $\alpha(q-2\ell)\pgq 0$ }\\
  (a\oa)^Ne_0\ot e_{1}(\oa a)^{N-1} & \mbox{if $q-2\ell=1$ and $\alpha<0$}\\
  (a\oa)^Ne_0\ot e_{-1}(a\oa)^{N-1} & \mbox{if $q-2\ell=-1$ and $\alpha>0$}\\
  0 & \mbox{otherwise}
  \end{cases}$\\for all $0\ppq\ell\ppq q$.
\item For each $j=0, 1, \ldots , m-2$ and each $s=1, \ldots , N-1$:\\
${\mathcal L}^q{F}_{n,j,s}:
  \begin{cases}
  e_j\ot e_{j+q-2\ell}\mapsto (a_j\oa_j)^se_j\ot e_{j+q-2\ell}\\
  e_{j+1}\ot e_{j+1+q-2\ell}\mapsto (-1)^{\frac{n}{2}}(\oa_ja_j)^s e_{j+1}\ot e_{j+1+q-2\ell}
  \end{cases}$\\for all $0\ppq\ell\ppq q$.
\item For each $s=1, \ldots , N-1$:\\
${\mathcal L}^q{F}_{n,m-1,s}:
  \begin{cases}
  e_{m-1}\ot e_{m-1+q-2\ell}\mapsto (a_{m-1}\oa_{m-1})^se_{m-1}\ot e_{m-1+q-2\ell}\\
  e_0\ot e_{q-2\ell}\mapsto (-1)^{\frac{n}{2}}(\oa_{m-1}a_{m-1})^s e_0\ot e_{q-2\ell}
  \end{cases}$\\for all $0\ppq\ell\ppq q$.
\item In the case $\car\kk\mid N$, for each
$j=1, \ldots , m-1$:\\
${\mathcal L}^q{\theta}_{n,j}:e_j\ot e_{j+q-2\ell}\mapsto
(a_j\oa_j)^Ne_j\ot e_{j+q-2\ell}$
  for all $0\ppq\ell\ppq q$.
\end{enumerate}

\item For $n$ odd:
\begin{enumerate}[$\bullet$]
\item For $-p\ppq \gamma<0$ and, if $t=m-1$, $\gamma=-p-1$ also:\\
${\mathcal L}^q{\varphi}_{n,\gamma}(e_i\ot e_{i+\gamma
m+1+q-2\ell})=
  \begin{cases}
  0&\mbox{if $q-2\ell>1$}\\
  -(-1)^{\frac{n-1-\gamma m}{2}i}(a\oa)^{N-1} e_i\ot e_{i+1}a(\oa a)^{N-1}&\mbox{if $q-2\ell=1$}\\
  (-1)^{q+\frac{n-1-\gamma m}{2}i}  e_i\ot e_{i+q-2\ell}a(\oa a)^{N-1}&\mbox{if $q-2\ell\ppq 0$}
  \end{cases}$\\
  for all $i=0, 1, \ldots , m-1$ and all $0\ppq\ell\ppq q$.
\item For $0\ppq \gamma \ppq p:$
\begin{enumerate}[$\diamond$]
      \item If $\gamma>0$ then ${\mathcal L}^q{\varphi}_{n,\gamma}(e_i\ot e_{i+\gamma m+1+q-2\ell})=\\
      \begin{cases}
      (-1)^{q+\frac{n-1-\gamma m}{2}i} e_i\ot e_{i+q-2\ell}a &\mbox{if $q-2\ell\pgq 0$}\\
      -(-1)^{\frac{n-1-\gamma m}{2}i}\left[\dsty\sum_{k=0}^{N-1}(\oa a)^ke_i\ot e_{i-1}(a\oa)^{N-k-1}a+\right.\\
      \hspace*{3cm}\left.(-1)^{\frac{q+1}{2}}\dsty\sum_{k=0}^{N-2}(\oa a)^k\oa e_{\rep{i-1}}\ot e_{\rep{i-1}+1} (\oa a)^{N-k-1}\right] &\mbox{if $q-2\ell=-1$}\\
      (-1)^{\frac{n-1-\gamma m}{2}i}(\oa a)^{N-1}e_i\ot e_{i-2}(a\oa)^{N-1}a &\mbox{if $q-2\ell=-2$}\\
      0&\mbox{if $q-2\ell<-2$}
      \end{cases}$\\for all $i=0, 1, \ldots , m-1$ and all $0\ppq\ell\ppq q$.

      \item If $\gamma=0$ then ${\mathcal L}^q{\varphi}_{n,0}(e_i\ot e_{i+1+q-2\ell})=\\
      \begin{cases}
      (-1)^{q+\frac{n-1}{2}i} e_i\ot e_{i+q-2\ell}a   &\mbox{if $q-2\ell=q$}\\
      (-1)^{q+\frac{n-1}{2}i} \left[e_i\ot e_{i+q-2\ell}a -(-1)^q(N-1)e_i\ot e_{i+q-2\ell+2}\oa(a\oa)^{N-1} \right] &\mbox{if $0\ppq q-2\ell<q$} \\
      (-1)^{q+\frac{n-1}{2}i} N e_i\ot e_{i+q-2\ell}a(\oa a)^{N-1}   &\mbox{if $q-2\ell<-1$} \\
      \!\!\!\begin{array}{l}
      -(-1)^{\frac{n-1}{2}i}
           \left\{\dsty\sum_{k=1}^{N-1}\sum_{v=0}^{k-1}
               \left[(a\oa)^ve_i\ot e_{i+1}\oa(a\oa)^{N-v-1}+ \right.\right.\\
               \hspace*{3.5cm}(-1)^{\frac{q+1}{2}}(a\oa)^vae_{\rep{i+1}}\ot e_{\rep{i+1}-1}(a\oa)^{N-v-1}+\\
               \hspace*{4cm}(-1)^{\frac{q+1}{2}}(\oa a)^v\oa e_{\rep{i-1}}\ot e_{\rep{i-1}+1}(\oa a)^{N-v-1}+\\
               \left.\hspace*{4.5cm}(\oa a)^ve_i\ot e_{i-1}a(\oa a)^{N-v-1}\right]\\
      \left.\hspace*{3cm}+\dsty\sum_{k=0}^{N-1}(\oa a)^ke_i\ot e_{i-1}a(\oa a)^{N-k-1}\right\}
      \end{array} &\mbox{if $q-2\ell=-1$}
      \end{cases}$\\for all $i=0, 1, \ldots , m-1$ and all $0\ppq\ell\ppq q$.
\end{enumerate}

\item For $0<\beta\ppq p$ and, if $t=m-1,$ $\beta=p+1$ also:\\
${\mathcal L}^q{\psi}_{n,\beta}(e_i\ot e_{i+\beta m-1+q-2\ell})=
  \begin{cases}
  0&\mbox{if $q-2\ell<-1$}\\
  (-1)^{\frac{n-1-\beta m}{2}i}(\oa a)^{N-1}e_i\ot e_{i-1}\oa(a\oa)^{N-1}&\mbox{if $q-2\ell=-1$}\\
  (-1)^{\frac{n-1-\beta m}{2}i} e_i\ot
  e_{i+q-2\ell}\oa(a\oa)^{N-1}&\mbox{if $q-2\ell>-1$}
  \end{cases}$\\for all $i=0, 1, \ldots , m-1$ and all $0\ppq\ell\ppq q$.

\item For $-p\ppq \beta\ppq 0$:
\begin{enumerate}[$\diamond$]
      \item If $\beta<0,$ then ${\mathcal L}^q{\psi}_{n,\beta}(e_i\ot e_{i+\beta m-1+q-2\ell})=\\
      \begin{cases}
      (-1)^{\frac{n-1-\beta m}{2}i}  e_i\ot e_{i+q-2\ell}\oa     &\mbox{if $q-2\ell<1$}\\
      (-1)^{\frac{n-1-\beta m}{2}i}  \left[\dsty\sum_{k=0}^{N-1}(a\oa)^k e_i\ot e_{i+1}\oa (a\oa)^{N-k-1}+\right.\\
      \hspace*{3cm}\left.(-1)^{\frac{q+1}{2}}\sum_{k=0}^{N-2}(a\oa)^kae_{\rep{i+1}}\ot e_{\rep{i+1}-1}(a\oa)^{N-k-1}\right]     &\mbox{if $q-2\ell=1$}\\
      (-1)^{\frac{n-1-\beta m}{2}i} (a\oa)^{N-1}e_i\ot e_{i+2}(\oa a)^{N-1}\oa      &\mbox{if $q-2\ell=2$}\\
      0 &\mbox{if $q-2\ell>2$}
      \end{cases}$\\for all $i=0, 1, \ldots , m-1$ and all $0\ppq\ell\ppq q$.

      \item If $\beta=0$ then ${\mathcal L}^q{\psi}_{n,0}(e_i\ot e_{i-1+q-2\ell})=\\
      \begin{cases}
      (-1)^{\frac{n-1}{2}i} e_i\ot e_{i+q-2\ell}\oa   &\mbox{if $q-2\ell=-q$}\\
      (-1)^{\frac{n-1}{2}i} \left[e_i\ot e_{i+q-2\ell}\oa -(-1)^q(N-1)e_i\ot e_{i+q-2\ell-2}a(\oa a)^{N-1} \right] &\mbox{if $-q< q-2\ell\ppq 0$} \\
      (-1)^{\frac{n-1}{2}i} N e_i\ot e_{i+q-2\ell}\oa(a \oa)^{N-1}   &\mbox{if $q-2\ell>1$} \\
      \begin{array}{l}
      \!\!\!(-1)^{\frac{n-1}{2}i}
          \left\{\dsty\sum_{k=1}^{N-1}\sum_{v=0}^{k-1}\left[(\oa a)^ve_i\ot e_{i-1}a(\oa a)^{N-v-1}+\right.\right.\\
          \hspace*{3cm}(-1)^{\frac{q+1}{2}}(\oa a)^v\oa e_{\rep{i-1}}\ot e_{\rep{i-1}+1}(\oa a)^{N-v-1}\\
          \hspace*{3.5cm}+(-1)^{\frac{q+1}{2}}(a \oa)^va e_{\rep{i+1}}\ot e_{\rep{i+1}-1}a \oa)^{N-v-1}+\\
          \hspace*{4cm}\left.(a \oa)^ve_i\ot e_{i+1}\oa(a \oa)^{N-v-1}\right]\\
      \left.\hspace*{2cm}+\dsty\sum_{k=0}^{N-1}(a \oa)^ke_i\ot e_{i+1}\oa(a \oa)^{N-k-1}\right\}
      \end{array} &\mbox{if $q-2\ell=1$}
      \end{cases}$\\for all $i=0, 1, \ldots , m-1$ and all $0\ppq\ell\ppq q$.
\end{enumerate}

\item For each $j=0, 1, \ldots , m-1$ and each $s = 1, \ldots, N-1$:\\
${\mathcal L}^q{E}_{n,j,s}=\\
  \begin{cases}
  e_j\ot e_{j+1}\mapsto  (a\oa)^s e_j\ot e_j a &\mbox{if $q$ is even}\\
  e_j\ot e_j\mapsto \dsty\sum_{k=0}^{N-s}\sum_{v=k}^{N-s}\left[(a\oa)^{v+s} e_j\ot e_{j+1}\oa(a\oa)^{N-v-1}+\right.\\
       \hspace*{3cm}\left.(-1)^{\frac{q+1}{2}}(a\oa)^{v+s}ae_{\rep{j+1}}\ot e_{\rep{j+1}-1}(a\oa)^{N-v-1}\right] &\mbox{if $q$ is odd}\\
 e_{\rep{j+1}}\ot e_{\rep{j+1}}\mapsto
  (-1)^{\frac{n-1}{2}}\dsty\sum_{k=0}^{N-s-1}\sum_{v=k}^{N-s}
  \left[(-1)^{\frac{q+1}{2}}(\oa a)^{v+s}\oa e_j\ot e_{j+1}(\oa a)^{N-v-1}\right.\\
  \hspace*{5.5cm}\left.+(\oa a)^{v+s+1}e_{\rep{j+1}}\ot e_{\rep{j+1}-1} a(\oa a)^{N-v-2}\right]
  & \mbox{if $q$ is odd}
  \end{cases}$

\item In the case $\car\kk\mid N$, for each
$j=1, \ldots , m-1$:\\
${\mathcal L}^q\omega_{n,j}=\\
  \begin{cases}
  e_j\ot e_j \mapsto \dsty\sum_{k=1}^{N-1}\sum_{v=0}^{k-1}
     \left[(-1)^{\frac{q-1}{2}}(a\oa)^vae_{j+1}\ot
     e_j(a\oa)^{N-v-1}-\right.\\
     \hspace*{3.5cm}\left.(a\oa)^ve_j\ot e_{j+1}\oa(a\oa)^{N-v-1}\right] \mbox{if $q$ is odd}\\
  e_{\rep{j+1}}\ot e_{\rep{j+1}}\mapsto (-1)^{\frac{n+1}{2}}\dsty\sum_{k=1}^{N-1}\sum_{v=0}^{k-1}
     \left[(\oa a)^{v+1}e_{\rep{j+1}}\ot e_{\rep{j+1}-1} a(\oa a)^{N-v-2}+\right.\\
     \hspace*{4.7cm}\left.(-1)^{\frac{q+1}{2}}(\oa a)^v\oa e_j\ot
       e_{j+1}(\oa a)^{N-v-1}\right]  \mbox{ if $q$ is odd}\\
  e_{\rep{j-v}}\ot e_{\rep{j-v}+2v}\mapsto (-1)^{\frac{n-1}{2}v}
     \left[(-1)^{v+\frac{q+1}{2}}ae_{\rep{j-v+1}}\ot e_{\rep{j-v+1}+2v-1}+\right.\\
     \hspace*{3.8cm}\left.e_{\rep{j-v}}\ot e_{\rep{j-v}+2v+1}\oa (a\oa)^{N-1}\right] \mbox{ for $1\ppq v\ppq \frac{q+1}{2}$ and $q$ odd}\\
  e_{\rep{j-v}}\ot e_{\rep{j-v}+2v+1}\mapsto (-1)^{\frac{n-1}{2}v}
     \left[e_{\rep{j-v}}\ot e_{\rep{j-v}+2v}a-\right.\\
     \hspace*{3.8cm}\left.(-1)^{v+\frac{q}{2}}\oa(a\oa)^{N-1}e_{\rep{j-v-1}}\ot e_{\rep{j-v-1}+2v+2}\right] \mbox{ for $0\ppq v\ppq \frac{q}{2}$
     and $q$ even.}
  \end{cases}$
\end{enumerate}
\end{enumerate}
\eprop

\bigskip

\subsection{The Hochschild cohomology ring $\HH^*(\L_1)$ for $N=1$, $m\pgq 3$, $m$ even}

\bt For $N=1$, $m\pgq 3$ and $m$ even, $\HH^*(\L_1)$ is a finitely
generated algebra with generators:\\
$$\begin{array}{ll}
1, a_i\oa_i \mbox{ for $i=0, 1, \ldots , m-1$} & \mbox{in degree 0},\\
\varphi_{1,0}, \psi_{1,0} & \mbox{in degree 1},\\
\chi_{2,0} & \mbox{in degree 2},\\
\varphi_{m-1,-1}, \psi_{m-1,1} & \mbox{in degree $m-1$},\\
\chi_{m,1}, \chi_{m,-1} & \mbox{in degree $m$}.
\end{array}$$
\et

It is known from \cite[Proposition 4.4]{SS} that those generators in
$\HH^*(\L_1)$ whose image is in $\rrad$ are nilpotent, and it may be
seen that the remaining generators of $\HH^*(\L_1)$, namely $1,
\chi_{2,0}, \chi_{m,1}$ and $\chi_{m,-1}$, are not nilpotent.
Moreover $\chi_{2,0}^m = \chi_{m,1}\chi_{m,-1}$. Thus we have the
following corollary.

\bc For $N=1$, $m\pgq 3$ and $m$ even,
$$\HH^*(\la)/\mathcal{N}\cong \kk[\chi_{2,0},\chi_{m,1},\chi_{m,-1}]/(\chi_{2,0}^m - \chi_{m,1}\chi_{m,-1})$$
and hence $\HH^*(\L_1)/\N$ is a commutative finitely generated
algebra of Krull dimension 2. \ec

\bigskip

\subsection{The Hochschild cohomology ring $\HH^*(\L_N)$ for $N>1$, $m\pgq 3$, $m$ even}

In this section we describe the main arguments behind the
presentation of $\HH^*(\L_N)$ by generators and relations. The first
lemma indicates the method required to show that $\HH^*(\L_N)$ is a finitely
generated algebra whilst the second lemma considers the relations
between the generators.

\bl\label{lem:gens} Suppose $N>1$, $m\pgq 3$ and $m$ even. Then
\begin{enumerate}
\item $\chi_{2,0}\chi_{2,0} = \chi_{4,0}$;
\item $\pi_{2,0} = (a_0\bar{a}_0)^N\chi_{2,0}$.
\end{enumerate}
\el

{\it Proof } (1) By definition, $\chi_{2,0}\chi_{2,0} =
\chi_{2,0}\cdot{\mathcal L}^2\chi_{2,0}$. Now, from Proposition
\ref{prop:liftings}, ${\mathcal L}^2\chi_{2,0}(e_i\otimes
e_{i+2-2\ell}) = (-1)^ie_i\otimes e_{i+2-2\ell}$ for $i = 0, \ldots
, m-1$ and $0 \ppq\ell\ppq 2$, and $\chi_{2,0}(e_i\otimes e_{i}) =
(-1)^ie_i$ for $i = 0, \ldots , m-1$. Thus
$$\chi_{2,0}\cdot{\mathcal L}^2\chi_{2,0}(e_i\otimes e_{i+2-2\ell})
= \left \{
\begin{array}{ll}(-1)^i(-1)^ie_i = e_i & \mbox{if $\ell = 1$,}\\
0 & \mbox{otherwise,}\end{array} \right. $$ for all $i = 0, \ldots ,
m-1$. Hence $\chi_{2,0}\chi_{2,0} = \chi_{4,0}$.

(2) We know that $\chi_{2,0}(e_i\otimes e_{i}) = (-1)^ie_i$ for $i =
0, \ldots , m-1$, $\pi_{2,0} : e_0\otimes e_0 \mapsto
(a_0\bar{a}_0)^N$, and $(a_0\bar{a}_0)^N \in \HH^0(\L_N)$. Then
$$((a_0\bar{a}_0)^N \chi_{2,0})(e_i\otimes e_{i}) = \chi_{2,0}((a_0\bar{a}_0)^Ne_i\otimes
e_{i}) =
\left \{ \begin{array}{ll} (-1)^0e_0(a_0\bar{a}_0)^N & \mbox{if $i=0$,}\\
0 & \mbox{if $i = 1, \ldots , m-1$.}\end{array} \right. $$ Thus
$\pi_{2,0} = (a_0\bar{a}_0)^N\chi_{2,0}$. \qed

\bigskip

We remark that, for $n$ even, $2 \ppq n < m$ and writing $n=pm+t$ as
usual, we necessarily have $p=0$, and thus, from Proposition
\ref{prop:basisHHn}, we have $\{\chi_{n,\alpha} \mid -p\ppq \alpha
\ppq p\} = \{\chi_{n,0}\}$. Inductively, Lemma \ref{lem:gens} (1)
gives, for $n=2k$ even, $2 \ppq n \ppq m-2$, that $\chi_{n,0} =
(\chi_{2,0})^k$.

\bigskip

\bl\label{lem:rels} Suppose $N>1$, $m\pgq 3$ and $m$ even. Then
\begin{enumerate}
\item $\chi_{m,1}\chi_{m,-1} = 0$;
\item $\varphi_{1,0}^2 =0$;
\item $\varphi_{1,0}\psi_{1,0} = Nm(a_0\oa_0)^N\chi_{2,0}$;
\item for $\car{\kk}\mid N$, $S = \sum_{k=1}^Nk$ and $\varepsilon_i = (a_i\oa_i)^N$
for $i=0, 1, \ldots , m-1$, we have
$$\varphi_{1,0}\omega_{1,i} =
\begin{cases}
S\chi_{2,0}(\varepsilon_i+\varepsilon_{i+1}) \mbox{ if $\car{\kk}=2$ and $i\neq m-1$}\\
S\chi_{2,0}(\varepsilon_{m-1}+\varepsilon_{0}) \mbox{ if $\car{\kk}=2$ and $i= m-1$}\\
0\mbox{ if $\car{\kk}\neq 2$.}
\end{cases}$$
\end{enumerate}
\el

{\it Proof } (1) We have $\chi_{m,1}\chi_{m,-1} =
\chi_{m,1}\cdot{\mathcal L}^m\chi_{m,-1}$. Since $m$ is even,
Proposition \ref{prop:liftings} gives
$${\mathcal L}^m\chi_{m,-1}(e_i\otimes e_{i+m-2\ell-m}) =
\begin{cases}
  e_i\ot e_{i+m-2\ell} & \mbox{if $m-2\ell\ppq 0$}\\
  0 & \mbox{if $m-2\ell > 0$ and $\abs{m-2\ell}>2$}\\
  (a\oa)^{N-1}e_i\ot e_{i+2}(\oa a)^{N-1} & \mbox{if $m-2\ell=2$}\\
  \end{cases}$$
for all $i=0, 1, \ldots , m-1$ and all $0\ppq\ell\ppq m$. Also
$\chi_{m,1}(e_i\otimes e_{i+m}) = e_i$ for $i=0, 1, \ldots , m-1$.
Hence $\chi_{m,1}\cdot{\mathcal L}^m\chi_{m,-1}(e_i\otimes
e_{i+m-2\ell-m}) = 0$ and $\chi_{m,1}\chi_{m,-1} = 0$.

(2) If $\car{\kk} \neq 2$, then $\varphi_{1,0}^2 =0$ by graded
commutativity of $\HH^*(\L_N)$ so we may assume $\car{\kk} = 2$.
Thus we have $\varphi_{1,0}(e_i \otimes e_{i+1}) = a_i$ for $i=0, 1,
\ldots , m-1$, and, from Proposition \ref{prop:liftings}, ${\mathcal
L}^1\varphi_{1,0}(e_i \otimes e_{i+1+1-2\ell}) =$
$$\begin{cases}
      e_i\ot e_{i+1}a_{i+1} & \mbox{if $\ell=0$,}\\
      \!\!\!\begin{array}{l}
      -\left\{\dsty\sum_{k=1}^{N-1}\sum_{v=0}^{k-1}
               \left[(a\oa)^ve_i\ot e_{i+1}\oa(a\oa)^{N-v-1} -
               (a\oa)^vae_{\rep{i+1}}\ot e_{\rep{i+1}-1}(a\oa)^{N-v-1}\right.\right.\\
               \left.\hspace*{2cm}-(\oa a)^v\oa e_{\rep{i-1}}\ot e_{\rep{i-1}+1}(\oa a)^{N-v-1} +
               (\oa a)^ve_i\ot e_{i-1}a(\oa a)^{N-v-1}\right]\\
      \left.\hspace*{3cm}+\dsty\sum_{k=0}^{N-1}(\oa a)^ke_i\ot e_{i-1}a(\oa a)^{N-k-1}\right\}
      \end{array} & \mbox{if $\ell=1$,}
      \end{cases}$$
for all $i=0, 1, \ldots , m-1$. Hence $$\varphi_{1,0}\cdot{\mathcal
L}^1\varphi_{1,0}(e_i \otimes e_{i+1+1-2\ell})=
\begin{cases}
      a_ia_{i+1} & \mbox{if $\ell=0$,}\\
      \!\!\!\begin{array}{l}
      -\dsty\sum_{k=1}^{N-1}\sum_{v=0}^{k-1}
               \left[(a\oa)^va_i\oa(a\oa)^{N-v-1} + (\oa a)^v\oa a_{i-1}(\oa a)^{N-v-1}
               \right]
      \end{array} & \mbox{if $\ell=1$.}
      \end{cases}$$
But $a_ia_{i+1} = 0$ and $(a_i\oa_i)^N + (\oa_{i-1}a_{i-1})^N = 0$
since $\car{\kk} = 2$. Thus $\varphi_{1,0}\cdot{\mathcal
L}^1\varphi_{1,0}(e_i \otimes e_{i+1+1-2\ell})= 0$ and so
$\varphi_{1,0}^2=0$.

(3) We have $\varphi_{1,0}(e_i \otimes e_{i+1}) = a_i$ for $i=0, 1,
\ldots , m-1$, and, from Proposition \ref{prop:liftings}, ${\mathcal
L}^1{\psi}_{1,0}(e_i\ot e_{i-1+1-2\ell})=$
$$\begin{cases}
      e_i\ot e_{i-1}\oa_{i-2} & \mbox{if $\ell=1$,}\\
      \begin{array}{l}
      \!\!\!\left\{\dsty\sum_{k=1}^{N-1}\sum_{v=0}^{k-1}\left[(\oa a)^ve_i\ot e_{i-1}a(\oa a)^{N-v-1} -
      (\oa a)^v\oa e_{\rep{i-1}}\ot e_{\rep{i-1}+1}(\oa a)^{N-v-1}\right.\right.\\
      \left.\hspace*{1.3cm} - (a \oa)^va e_{\rep{i+1}}\ot e_{\rep{i+1}-1}a \oa)^{N-v-1} +
      (a \oa)^ve_i\ot e_{i+1}\oa(a \oa)^{N-v-1}\right]\\
      \left.\hspace*{2cm}+\dsty\sum_{k=0}^{N-1}(a \oa)^ke_i\ot e_{i+1}\oa(a \oa)^{N-k-1}\right\}
      \end{array} &\mbox{if $\ell=0$,}
      \end{cases}$$
for $i=0, 1, \ldots , m-1$. Thus $\varphi_{1,0}\cdot{\mathcal
L}^1{\psi}_{1,0}(e_i\ot e_{i-1+1-2\ell})=0$ if $\ell = 1$. On the
other hand, if $\ell=0$, then $\varphi_{1,0}\cdot{\mathcal
L}^1{\psi}_{1,0}(e_i\ot e_i) =
\dsty\sum_{k=1}^{N-1}\sum_{v=0}^{k-1}\left[ -(\oa a)^v\oa
a_{i-1}(\oa a)^{N-v-1} + (a \oa)^va_i\oa(a \oa)^{N-v-1}\right] +
\dsty\sum_{k=0}^{N-1}(a \oa)^ka_i\oa(a\oa)^{N-k-1} =
\dsty\sum_{k=0}^{N-1}(a_i\oa_i)^N = N(a_i\oa_i)^N$.

Now define $g \in \Hom(P^1, \L_N)$ by $e_i\otimes e_{i+1} \mapsto
(m-1-i)a_i$ for $i = 0, 1, \ldots , m-2$ (with our usual convention
that all other summands are mapped to zero). Then $\partial^2g \in
\Hom(P^2, \L_N)$ and, for $i = 1, \ldots , m-2$, we have
$$\begin{array}{ll}
\partial^2g(e_i\otimes e_i) & = \sum_{k=0}^{N-1}[-(\oa a)^k\oa g(e_{i-1}\ot e_i)(\oa a)^{N-k-1}
+ (a\oa)^kg(e_i\ot e_{i+1})\oa(a\oa)^{N-k-1}]\\
& = \sum_{k=0}^{N-1}[-(\oa a)^k\oa (m-i)a_{i-1}(\oa a)^{N-k-1}
+ (a\oa)^k(m-1-i)a_i\oa(a\oa)^{N-k-1}]\\
& = \sum_{k=0}^{N-1}[-(m-i)(\oa_{i-1}a_{i-1})^N + (m-1-i)(a_i\oa_i)^N]\\
& = -N(a_i\oa_i)^N.
\end{array}$$
In a similar way, we also have $\partial^2g(e_{m-1}\otimes e_{m-1})
= -N(a_{m-1}\oa_{m-1})^N$. For $i = 0$ we have,
$$\begin{array}{ll}
\partial^2g(e_0\otimes e_0) & =
\sum_{k=0}^{N-1}[-(\oa a)^k\oa g(e_{m-1}\ot e_m)(\oa a)^{N-k-1}
+ (a\oa)^kg(e_0\ot e_1)\oa(a\oa)^{N-k-1}]\\
& = \sum_{k=0}^{N-1}(a\oa)^k(m-1)a_0\oa(a\oa)^{N-k-1}\\
& = N(m-1)(a_0\oa_0)^N.
\end{array}$$
Moreover, $\partial^2g$ is zero on all other summands. Hence
$$\partial^2g(e_i\otimes e_i) =
\left \{ \begin{array}{ll}N(m-1)(a_0\oa_0)^N & \mbox{if $i=0$},\\
-N(a_i\oa_i)^N & \mbox{if $i = 1, \ldots , m-1$}.
\end{array}\right.$$
So $\varphi_{1,0}\cdot {\mathcal L}^1\psi_{1,0} + \partial^2g$ is
the map given by $e_0 \otimes e_0 \mapsto Nm(a_0\oa_0)^N$ and thus
$\varphi_{1,0}\cdot {\mathcal L}^1\psi_{1,0} + \partial^2g =
Nm\pi_{2,0}$. Using Lemma \ref{lem:gens}(2), we have that the
cocycles $\varphi_{1,0}\psi_{1,0}$ and $Nm(a_0\oa_0)^N\chi_{2,0}$
represent the same elements in Hochschild cohomology so that
$\varphi_{1,0}\psi_{1,0} = Nm(a_0\oa_0)^N\chi_{2,0}$ as required.

(4) Finally, suppose $\car{\kk}\mid N$, $S = \sum_{k=1}^Nk$ and let
$\varepsilon_i = (a_i\oa_i)^N$ for $i=0, 1, \ldots , m-1$. For ease
of notation we consider the product $\varphi_{1,0}\omega_{1,i}$
where $1 \ppq i \ppq m-2$; the cases $i=0$ and $i = m-1$ are
similar. We know $\varphi_{1,0}(e_j\otimes e_{j+1}) = a_j$ for $j =
0, 1, \ldots , m-1$ and, from Proposition \ref{prop:liftings}, we
have ${\mathcal L}^1\omega_{1,i}:
\begin{cases}
  e_i\ot e_i \mapsto \dsty\sum_{k=1}^{N-1}\sum_{v=0}^{k-1}
     \left[(a\oa)^vae_{i+1}\ot e_i(a\oa)^{N-v-1} - (a\oa)^ve_i\ot e_{i+1}\oa(a\oa)^{N-v-1}\right] \\
  e_{i+1}\ot e_{i+1}\mapsto -\dsty\sum_{k=1}^{N-1}\sum_{v=0}^{k-1}
     \left[(\oa a)^{v+1}e_{i+1}\ot e_i a(\oa a)^{N-v-2} - (\oa a)^v\oa e_i\ot e_{i+1}(\oa a)^{N-v-1}\right] \\
  e_{i-1}\ot e_{i+1}\mapsto ae_{i}\ot e_{i+1}.
\end{cases}$
Hence $\varphi_{1,0}\cdot{\mathcal L}^1\omega_{1,i}(e_i\ot e_i) =
\dsty\sum_{k=1}^{N-1}\sum_{v=0}^{k-1} -(a\oa)^va_i\oa(a\oa)^{N-v-1}
= \dsty\sum_{k=1}^{N-1}-k(a_i\oa_i)^N = -S(a_i\oa_i)^N$. And
$\varphi_{1,0}\cdot{\mathcal L}^1\omega_{1,i}(e_{i+1}\ot e_{i+1})=
(-1)\dsty\sum_{k=1}^{N-1}\sum_{v=0}^{k-1}(-1)(\oa a)^v\oa a_i(\oa
a)^{N-v-1} = \dsty\sum_{k=1}^{N-1}k(\oa_i a_i)^N = S(\oa_i a_i)^N$.
Finally, $\varphi_{1,0}\cdot{\mathcal L}^1\omega_{1,i}(e_{i-1}\ot
e_{i+1}) = 0$. Hence $\varphi_{1,0}\cdot{\mathcal L}^1\omega_{1,i}:
\begin{cases}
e_i\otimes e_i \mapsto -S(a_i\oa_i)^N\\
e_{i+1}\ot e_{i+1} \mapsto S(a_{i+1} \oa_{i+1})^N.
\end{cases}$
Now, $\chi_{2,0}(\varepsilon_i + \varepsilon_{i+1}):
\begin{cases}
e_i\ot e_i \mapsto (-1)^i(a_i\oa_i)^N\\
e_{i+1}\ot e_{i+1} \mapsto (-1)^{i+1}(a_{i+1} \oa_{i+1})^N.
\end{cases}$
If $\car{\kk} = 2$, then $\varphi_{1,0}\omega_{1,i} =
\varphi_{1,0}\cdot{\mathcal L}^1\omega_{1,i} =
S\chi_{2,0}(\varepsilon_i + \varepsilon_{i+1})$. However, if
$\car{\kk} \neq 2$ then, since $\car{\kk} \mid N$, we have $2S =
(N-1)N = 0$ so that $S=0$ and hence $\varphi_{1,0}\omega_{1,i} = 0$.
\qed

\bigskip

By considering all possible products of the basis elements of
$\HH^*(\L_N)$ from Proposition \ref{prop:basisHHn}, we are now able
to give the main result of this section. We do not explicitly list
the relations which are a direct consequence of the graded
commutativity of $\HH^*(\L_N)$. Note that we can check that we do
indeed have all the relations by comparing $\dim\HH^n(\L_N)$ from
Proposition \ref{prop:dim} with the dimension in degree $n$ of the
algebra given in Theorem \ref{thm:thm!} below.

\bigskip

\bt\label{thm:thm!} Suppose $N>1$, $m\pgq 3$ and $m$ even.
\begin{enumerate}
\item $\HH^*(\L_N)$ is a finitely
generated algebra with generators:\\
$$\begin{array}{ll}
1, (a_i\oa_i)^N , [a_i\oa_i+\oa_ia_i] \mbox{ for $i=0, 1, \ldots , m-1$} & \mbox{in degree 0}\\
\varphi_{1,0}, \psi_{1,0} & \mbox{in degree 1}\\
\omega_{1,j} \mbox{ for $j=1, \ldots , m-1$ if
$\car\kk\mid N$} & \mbox{in degree 1}\\
\chi_{2,0} & \mbox{in degree 2}\\
\varphi_{m-1,-1}, \psi_{m-1,1} & \mbox{in degree $m-1$}\\
\chi_{m,1}, \chi_{m,-1} & \mbox{in degree $m$.}
\end{array}$$

\item Let $\varepsilon_i = (a_i\oa_i)^N$ for $i=0, 1, \ldots , m-1$ and $S = \sum_{k=1}^Nk$.
Then $\HH^*(\L_N)$ is a finitely
generated algebra over $\HH^0(\L_N)$ with generators:\\
$$\begin{array}{ll}
1 & \mbox{in degree 0}\\
\varphi_{1,0}, \psi_{1,0} & \mbox{in degree 1}\\
\omega_{1,j} \mbox{ for $j=1, \ldots , m-1$ if
$\car\kk\mid N$} & \mbox{in degree 1}\\
\chi_{2,0} & \mbox{in degree 2}\\
\varphi_{m-1,-1}, \psi_{m-1,1} & \mbox{in degree $m-1$}\\
\chi_{m,1}, \chi_{m,-1} & \mbox{in degree $m$,}
\end{array}$$
and relations
$$
\begin{array}{lll}
\varphi_{1,0}^2=0 \ & \ \varphi_{1,0}\psi_{1,0}=Nm\varepsilon_0\chi_{2,0} \ & \ \varphi_{1,0}\varphi_{m-1,-1}=0  \\
\psi_{1,0}^2=0  \ & \ \psi_{1,0}\psi_{m-1,1}=0  \ & \ \varphi_{m-1,-1}^2=0   \\
\varphi_{m-1,-1}\psi_{m-1,1}=0 \ & \
\varphi_{m-1,-1}\chi_{m,1}=0 \
&\ \psi_{m-1,1}^2=0  \\
\psi_{m-1,1}\chi_{m,-1}=0 \ & \ \chi_{m,1}\chi_{m,-1}=0 \ & \
\psi_{1,0}\chi_{m,1} = N\chi_{2,0}\psi_{m-1,1} \\
\varphi_{1,0}\chi_{m,-1} = N\chi_{2,0}\varphi_{m-1,-1} \ & \
\omega_{1,j}\psi_{m-1,1} =
(-1)^{\frac{m}{2}j}\chi_{m,1}\varepsilon_0 \ & \ \omega_{1,j}\chi_{m,\pm1} = 0 \\
\omega_{1,i}\omega_{1,j} = 0\mbox{ if $i\neq j$} &
\multicolumn{2}{l}{\ \varphi_{1,0}\psi_{m-1,1} = m\chi_{m,1}\varepsilon_0 = \psi_{1,0}\varphi_{m-1,-1}}\\
\multicolumn{3}{l}{\varphi_{1,0}\omega_{1,i} =
\psi_{1,0}\omega_{1,i} = \omega_{1,i}^2 =
\begin{cases}
S\chi_{2,0}(\varepsilon_i+\varepsilon_{i+1}) \mbox{ if $\car{\kk}=2$ and $i\neq m-1$}\\
S\chi_{2,0}(\varepsilon_{m-1}+\varepsilon_{0}) \mbox{ if $\car{\kk}=2$ and $i= m-1$}\\
0\mbox{ if $\car{\kk}\neq 2$.}
\end{cases}}
\end{array}$$

\item Let $\varepsilon_i = (a_i\oa_i)^N$ and $f_i = [a_i\oa_i+\oa_ia_i]$
for $i=0, 1, \ldots , m-1$. Then the action of $\HH^0(\L_N)$ on the
generators of $\HH^n(\L_N)$ is described as follows, for all $i, j =
0, 1, \ldots , m-1$.
$$\begin{array}{lll} \varepsilon_ig = 0 \ \forall g\neq\chi_{n,\alpha}
& \ \varepsilon_i\chi_{2,0} = 0 \mbox{ if $i\neq0$ and
$\car{\kk}\nmid N$} & \ \varepsilon_i\chi_{m,\pm 1} = (-1)^i\varepsilon_0\chi_{m,\pm1}\\
f_if_j = 0 \mbox{ if $i\neq j$} \ & \ f_i\varphi_{1,0} =
f_i\psi_{1,0} & \ f_i^N = \varepsilon_i+\varepsilon_{i+1}\\
f_i\varphi_{m-1,-1} = 0 & \ f_i\psi_{m-1,1} = 0 & \ f_i\chi_{m,\pm1} = 0\\
f_i\omega_{1,j} = 0 \mbox{ for $i\neq j$}\ & \ f_i\omega_{1,i} =
f_i\varphi_{1,0}. &
\end{array}$$

\end{enumerate}
\et

\bc For $N>1$, $m\pgq 3$ and $m$ even,
$$\HH^*(\L_N)/\mathcal{N}\cong \kk[\chi_{2,0},\chi_{m,1},\chi_{m,-1}]/(\chi_{m,1}\chi_{m,-1})$$
and hence $\HH^*(\L_1)/\N$ is a commutative finitely generated
algebra of Krull dimension 2. \ec

\bigskip

\section{The Hochschild cohomology ring $\HH^*(\L_N)$ for $m \pgq 3$, $m$ odd.}

In this section, $N \pgq 1$, $m \pgq 3$ and $m$ is odd. Most of the
details are similar to those in the previous section and are left to
the reader. We consider separately the cases $\car\kk \neq 2$ and
$\car\kk = 2$.

\bigskip

\subsection{Basis of $\HH^n(\L_N)$ for $N \pgq 1$, $m\pgq 3$, $m$ odd and $\car\kk \neq 2$.}

\bprop Suppose that $N\pgq 1$, $m\pgq 3$, $m$ odd and $\car\kk \neq
2$. For $n\pgq 1$, the following elements define a basis of
$\HH^n(\L_N)$:
\begin{enumerate}
\item For $n$ even, $n \pgq 2$:
\begin{enumerate}
\item $\chi_{n,\delta}: \ e_i\ot e_{i+\delta m}\mapsto e_i$ for all $i=0, 1, \ldots , m-1$, with\\
$\delta= \begin{cases}
p-2\alpha-1,\ 0\ppq \alpha\ppq p-1 & \mbox{when $t$ is odd and $\alpha+\frac{m-t}{2}$ is odd},\\
p-2\alpha, \ 0\ppq \alpha \ppq p & \mbox{when $t$ is even and
$\alpha+\frac{t}{2}$ is even; }
\end{cases}$
\item $\pi_{n,\delta}:\ e_0\ot e_{\delta m}\mapsto (a_0\oa_0)^N$, with\\
$\delta = \begin{cases}
p-2\alpha-1,\ 0\ppq \alpha\ppq p-1 & \mbox{when $t$ is odd and $\alpha+\frac{m-t}{2}$ is even},\\
p-2\alpha, \ 0\ppq \alpha \ppq p & \mbox{when $t$ is even and
$\alpha+\frac{t}{2}$ is odd; }
\end{cases}$
\item For each $j=0, 1, \ldots , m-2$ and each $s = 1, \ldots , N-1$:\\ $F_{n,j,s}:
\begin{cases}
e_j\ot e_j\mapsto (a_j\oa_j)^s,\\
e_{j+1}\ot e_{j+1}\mapsto (-1)^\frac{n}{2}(\oa_ja_j)^s;
\end{cases}$
\item For each $s = 1, \ldots , N-1$:\\ $F_{n,m-1,s}:
\begin{cases}
e_{m-1}\ot e_{m-1}\mapsto (a_{m-1}\oa_{m-1})^s,\\
e_0\ot e_0\mapsto (-1)^\frac{n}{2}(\oa_{m-1}a_{m-1})^s;
\end{cases}$
\item For $t=m-1$ and $\sigma=-(p+1)$: $\varphi_{n,\sigma}:\ e_i\ot e_{i+\sigma m+1}\mapsto
(a\oa)^{N-1}a_i$ for all $i=0, 1, \ldots , m-1$;
\item For $t=m-1$ and $\tau=p+1$: $\psi_{n,\tau}:\ e_i\ot e_{i+\tau m-1}\mapsto
(\oa a)^{N-1}\oa_{i-1}$ for all $i=0, 1, \ldots , m-1$;
\item Additionally, if $\car\kk\mid N$, and for each
$\begin{cases}
j=0, 1, \ldots , m-1 & \mbox{if $\frac{n}{2}$ is even}\\
j=1, \ldots , m-1 & \mbox{if $\frac{n}{2}$ is odd:}
\end{cases}$\\
$\theta_{n,j}:e_j\ot e_j\mapsto (a_j\oa_j)^N$.
\end{enumerate}

\item For $n$ odd:
\begin{enumerate}
\item For $t=0$ and $\delta = \pm p$:
$\pi_{n,\delta}:\ e_0\ot e_{\delta m}\mapsto (a_0\oa_0)^N$;
\item $\varphi_{n,\sigma}: e_i\ot e_{i+\sigma m+1}\mapsto (a\oa)^{N-1}a_i$ for all $i=0, 1, \ldots , m-1$, with\\
$\begin{cases}
\sigma=p-2\gamma & \mbox{if $t$ is odd, $2\gamma>p$, $\gamma\ppq p$, $\gamma+\frac{t-1}{2}$ is even,}\\
\sigma=p-2\gamma-1 & \mbox{if $t$ is even, $t\neq m-1$, $\gamma\ppq
p-1$,
  $2\gamma>p-1$ and $\gamma+\frac{m-1}{2}+\frac{t}{2}$ even,}\\
\sigma=p-2\gamma-1 & \mbox{if $t=m-1,$ $ \gamma\ppq p,$
$2\gamma>p-1$, $\gamma$ even;}
\end{cases}$
\item $\varphi_{n,\sigma}: e_i\ot e_{i+\sigma m+1}\mapsto a_i$ for all $i=0, 1, \ldots , m-1$, with \\
$\begin{cases}
\sigma=p-2\gamma & \mbox{if $t$ is odd, $2\gamma\ppq p$, $\gamma\pgq 0$  and $\gamma+\frac{t-1}{2}$ even,}\\
\sigma=p-2\gamma-1 & \mbox{if $t$ is even, $ \gamma\pgq 0$, $2\gamma\ppq p-1$ and $\gamma+\frac{m-1}{2}+\frac{t}{2}$ even;}\\
\end{cases}$
\item $\psi_{n,\tau}: e_i\ot e_{i+\tau m-1}\mapsto (\oa a)^{N-1}\oa_{i-1}$ for all $i=0, 1, \ldots , m-1$, with\\
$\begin{cases}
\tau=p-2\beta & \mbox{if $t$ is odd, $2\beta<p$, $\beta\pgq 0$, and $\beta+\frac{t-1}{2}$ even,}\\
\tau=p-2\beta-1 & \mbox{if $t$ is even, $t\neq m-1$, $\beta\pgq 0$,
  $2\beta<p-1$ and $\beta+\frac{m-1}{2}+\frac{t}{2}$ even, }\\
\tau=p-2\beta-1 & \mbox{if $t=m-1$, $\beta\pgq -1,$ $2\beta<p-1,$ and $\beta$ even;}\\
\end{cases}$
\item $\psi_{n,\tau}: e_i\ot e_{i+\tau m-1}\mapsto \oa_{i-1}$ for all $i=0, 1, \ldots , m-1$,  with\\
$\begin{cases} \tau=p-2\beta & \mbox{if $t$ is odd, $2\beta\pgq p$,
  $\beta\ppq p$, $2\beta\pgq p$, and $\beta+\frac{t-1}{2}$ even,}\\
\tau=p-2\beta-1 & \mbox{if $t$ is even, $ \beta\ppq p-1$,
$2\beta\pgq p-1$,
and $\beta+\frac{m-1}{2}+\frac{t}{2}$ even;}\\
\end{cases}$
\item For each $j=0, 1, \ldots , m-1$ and each $s=1, \ldots , N-1$: $E_{n,j,s}:\ e_j\ot
e_{j+1}\mapsto (a\oa)^sa_j$;
\item Additionally, if $\car\kk\mid N$, and for each
$\begin{cases}
j=0, 1, \ldots , m-1 & \mbox{if $\frac{n-1}{2}$ is odd}\\
j=1, \ldots , m-1 & \mbox{if $\frac{n-1}{2}$ is even:}
\end{cases}$\\
$\omega_{n,j}:e_j\ot e_{j+1}\mapsto a_j$.
\end{enumerate}
\end{enumerate}
\eprop

\subsection{The Hochschild cohomology ring $\HH^*(\L_N)$ for $N \pgq 1$, $m\pgq 3$, $m$ odd and $\car\kk \neq 2$}

\bt For $N=1$, $m\pgq 3$, $m$ odd and $\car\kk\neq 2$, $\HH^*(\L_1)$
is a finitely generated algebra with generators:\\
$$\begin{array}{ll}
1, a_i\oa_i \mbox{ for $i=0, 1, \ldots , m-1$} & \mbox{in degree 0,}\\
\varphi_{1,0}, \psi_{1,0}  & \mbox{in degree 1,}\\
\pi_{2,0} & \mbox{in degree 2,}\\
\chi_{4,0} & \mbox{in degree 4,}\\
\varphi_{m-1,-1}, \psi_{m-1,1} & \mbox{in degree $m-1$,}\\
\pi_{m,1}, \pi_{m,-1} & \mbox{in degree $m$,}\\
\chi_{2m,2}, \chi_{2m,-2} & \mbox{in degree $2m$.}
\end{array}$$
\et

We note that if $\car(\kk)\nmid m$, then the generators $\pi_{2,0},
\pi_{m,1}$ and $\pi_{m,-1}$ are redundant.

\bc For $N=1$, $m\pgq 3$, $m$ odd and $\car\kk\neq 2$,
$$\HH^*(\L_1)/\mathcal{N}\cong \kk[\chi_{4,0},\chi_{2m,2},\chi_{2m,-2}]/(\chi_{4,0}^m - \chi_{2m,2}\chi_{2m,-2})$$
and hence $\HH^*(\L_1)/\N$ is a commutative finitely generated
algebra of Krull dimension 2. \ec

\bt For $N>1$, $m\pgq 3$, $m$ odd and $\car\kk\neq 2$, $\HH^*(\L_N)$
is a finitely
generated algebra with generators:\\
$$\begin{array}{ll}
1, (a_i\oa_i)^N, [a_i\oa_i+\oa_ia_i] \mbox{ for $i=0, 1, \ldots , m-1$} & \mbox{in degree 0,}\\
\varphi_{1,0}, \psi_{1,0} & \mbox{in degree 1,}\\
\omega_{1,j} \mbox{ for $j=1, \ldots , m-1$ if
$\car\kk\mid N$} & \mbox{in degree 1,}\\
\pi_{2,0}, F_{2,j,1} \mbox{ for $j=1, \ldots , m-1$} & \mbox{in degree 2}\\
\theta_{2,j} \mbox{ for $j=1, \ldots , m-1$ if
$\car\kk\mid N$} & \mbox{in degree 2,}\\
\omega_{3,j} \mbox{ for $j=1, \ldots , m-1$ if
$\car\kk\mid N$} & \mbox{in degree 3,}\\
\chi_{4,0} & \mbox{in degree 4,}\\
\varphi_{m-1,-1}, \psi_{m-1,1} & \mbox{in degree $m-1$,}\\
\pi_{m,1}, \pi_{m,-1} & \mbox{in degree $m$,}\\
\chi_{2m,2}, \chi_{2m,-2} & \mbox{in degree $2m$,}\\
\varphi_{2m+1,-2}, \psi_{2m+1,2} \mbox{ if
$\car\kk\mid N$} & \mbox{in degree $2m+1$.}\\
\end{array}$$
\et

\bc For $N>1$, $m\pgq 3$, $m$ odd and $\car\kk\neq 2$,
$$\HH^*(\L_N)/\mathcal{N}\cong \kk[\chi_{4,0},\chi_{2m,2},\chi_{2m,-2}]/(\chi_{2m,2}\chi_{2m,-2})$$
and hence $\HH^*(\L_1)/\N$ is a commutative finitely generated
algebra of Krull dimension 2. \ec

\bigskip

\subsection{Basis of $\HH^n(\L_N)$ for $N\pgq 1$, $m\pgq 3$, $m$ odd and $\car\kk = 2$.}

\bprop Suppose that $N\pgq 1$, $m\pgq 3$, $m$ odd and $\car\kk = 2$.
For $n \pgq 1$, the following elements define a basis of
$\HH^n(\L_N)$:
\begin{enumerate}
\item For all $n \pgq 1$:
\begin{enumerate}
\item $\chi_{n,\delta}: e_i\ot e_{i+\delta m}\mapsto e_i$ for all
$i=0, 1, \ldots , m-1$, with\\ $\delta=
\begin{cases}
p-2\alpha-1,\ 0\ppq \alpha\ppq p-1 & \mbox{when $t$ is odd},\\
p-2\alpha, \ 0\ppq \alpha \ppq p & \mbox{when $t$ is even;}
\end{cases}$
\item $\pi_{n,\delta}: e_0\ot e_{\delta m}\mapsto (a_0\oa_0)^N$, with $\delta=
\begin{cases}
p-2\alpha-1,\ 0\ppq \alpha\ppq p-1 & \mbox{when $t$ is odd},\\
p-2\alpha, \ 0\ppq \alpha \ppq p & \mbox{when $t$ is even;}
\end{cases}$
\item $\varphi_{n,\sigma}: e_i\ot e_{i+\sigma m+1}\mapsto
a_i$ for all $i = 0, 1, \dots , m-1$ with\\ $\sigma =
\begin{cases}
p-2\gamma & \mbox{if $t$ is odd, $2\gamma\ppq p$, and $\gamma\pgq 0$, }\\
p-2\gamma-1 & \mbox{if $t$ is even, $t\neq m-1,$ $\gamma\pgq 0$, and $2\gamma\ppq p-1$,}\\
p-2\gamma-1 & \mbox{if $t=m-1,$ $2\gamma\ppq p-1$, $\gamma\pgq 0$;}
\end{cases}$
\item $\varphi_{n,\sigma}: e_i\ot e_{i+\sigma m+1}\mapsto
(a\oa)^{N-1}a_i$ for all $i=0, 1, \ldots , m-1$ with\\ $\sigma =
\begin{cases}
p-2\gamma & \mbox{if $t$ is odd, $2\gamma>p$, and $\gamma\ppq p$, }\\
p-2\gamma-1 & \mbox{if $t$ is even, $t\neq m-1,$ $ \gamma\ppq p-1$, and $2\gamma>p-1$,}\\
p-2\gamma-1 & \mbox{if $t=m-1,$ $2\gamma>p-1,$ $\gamma\ppq p$;}
\end{cases}$
\item $\psi_{n,\tau}: e_i\ot e_{i+\tau m-1}\mapsto
\oa_{i-1}$ for all $i=0, 1, \ldots , m-1$, with:\\
$\tau=\begin{cases}
p-2\beta & \mbox{if $t$ is odd, $2\beta\pgq p$, and $\beta\ppq p$,}\\
p-2\beta-1 & \mbox{if $t$ is even, $t\neq m-1,$ $\beta\ppq p-1$, and $2\beta\pgq p-1$,}\\
p-2\beta-1 & \mbox{if $t=m-1,$ $2\beta\pgq p-1,$ $\beta\ppq p-1$;}
\end{cases}$
\item $\psi_{n,\tau}: e_i\ot e_{i+\tau m-1}\mapsto
(\oa a)^{N-1}\oa_{i-1}$ for all $i=0, 1, \ldots , m-1$, with:\\
$\tau = \begin{cases}
p-2\beta & \mbox{if $t$ is odd, $2\beta<p$, and $\beta\pgq 0$,}\\
p-2\beta-1 & \mbox{if $t$ is even, $t\neq m-1,$ $ \beta\pgq 0$, and $2\beta<p-1,$}\\
p-2\beta-1 & \mbox{if $t=m-1,$ $2\beta<p-1$, $\beta\pgq -1.$}
\end{cases}$
\end{enumerate}

\item For $n$ even, $n\pgq 2$:
\begin{enumerate}
\item For each $j = 0, 1, \ldots , m-2$ and each $s=1, \ldots , N-1$:\\ $F_{n,j,s}:
\begin{cases}
e_j\ot e_j\mapsto (a_j\oa_j)^s,\\
e_{j+1}\ot e_{j+1}\mapsto (\oa_ja_j)^s;
\end{cases}$
\item For each $s=1, \ldots , N-1$:\\ $F_{n,m-1,s}:
\begin{cases}
e_{m-1}\ot e_{m-1}\mapsto (a_{m-1}\oa_{m-1})^s,\\
e_0\ot e_0\mapsto (\oa_{m-1}a_{m-1})^s;
\end{cases}$
\item Additionally, if $\car\kk\mid N$ and for each $j=0, 1, \ldots ,
m-1$: $\theta_{n,j}: e_j\ot e_j\mapsto (a_j\oa_j)^N$.
\end{enumerate}

\item For $n$ odd:
\begin{enumerate}
\item For each $j=0, 1, \ldots , m-1$ and each $s = 1, \ldots, N-1$:
$E_{n,j,s}: e_j\ot e_{j+1}\mapsto (a\oa)^sa_j$;
\item For each $j=0, 1, \ldots , m-1$: $\omega_{n,j}: e_j\ot e_{j+1}\mapsto
a_j$.
\end{enumerate}
\end{enumerate}
\eprop

\bigskip

\subsection{The Hochschild cohomology ring $\HH^*(\L_1)$ for $N\pgq 1$, $m\pgq 3$, $m$ odd and $\car\kk = 2$}

\bt For $N=1$, $m\pgq 3$, $m$ odd and $\car\kk = 2$, $\HH^*(\L_1)$
is a finitely generated algebra with generators:\\
$$\begin{array}{ll}
1, a_i\oa_i \mbox{ for $i=0, 1, \ldots , m-1$} & \mbox{in degree 0,}\\
\varphi_{1,0}, \psi_{1,0} & \mbox{in degree 1,}\\
\chi_{2,0} & \mbox{in degree 2,}\\
\varphi_{m-1,-1}, \psi_{m-1,1} & \mbox{in degree $m-1$,}\\
\chi_{m,1}, \chi_{m,-1} & \mbox{in degree $m$.}
\end{array}$$
\et

\bc For $N=1$, $m\pgq 3$, $m$ odd and $\car\kk = 2$,
$$\HH^*(\L_1)/\mathcal{N}\cong \kk[\chi_{2,0},\chi_{m,1},\chi_{m,-1}]/(\chi_{2,0}^m - \chi_{m,1}\chi_{m,-1})$$
and hence $\HH^*(\L_1)/\N$ is a commutative finitely generated
algebra of Krull dimension 2. \ec

\bt For $N>1$, $m\pgq 3$, $m$ odd and $\car\kk = 2$, $\HH^*(\L_N)$
is a finitely generated algebra with generators:\\
$$\begin{array}{ll}
1, (a_i\oa_i)^N , [a_i\oa_i+\oa_ia_i] \mbox{ for $i=0, 1, \ldots , m-1$} & \mbox{in degree 0,}\\
\varphi_{1,0}, \psi_{1,0} & \mbox{in degree 1,}\\
\omega_{1,j} \mbox{ for $j=0, 1, \ldots , m-1$ if $\car\kk\mid N$} & \mbox{in degree 1,}\\
\chi_{2,0} & \mbox{in degree 2,}\\
\varphi_{m-1,-1}, \psi_{m-1,1} & \mbox{in degree $m-1$,}\\
\chi_{m,1}, \chi_{m,-1} & \mbox{in degree $m$.}
\end{array}$$
\et

\bc For $N>1$, $m\pgq 3$, $m$ odd and $\car\kk = 2$,
$$\HH^*(\lan)/\mathcal{N}\cong \kk[\chi_{2,0}, \chi_{m,1}, \chi_{m,-1}]/(\chi_{m,1}\chi_{m,-1})$$
and hence $\HH^*(\L_N)/\N$ is a commutative finitely generated
algebra of Krull dimension 2. \ec

\bigskip

\section{The Hochschild cohomology ring $\HH^*(\L_N)$ for $m=2$}

Recall, from Theorem \ref{thm:centre}, that $\dim\HH^0(\lan) =
2N+1$. We start with the case $N = 1$.

\subsection{The case $N=1$}

\bprop For $m=2$ and $N=1$,
$$\dim\HH^n(\la)=
\begin{cases}
3&\mbox{ if $n=0,$}\\
2(n+1)&\mbox{ if $n\neq 0.$}
\end{cases}$$
\eprop

\bprop For $m=2$ and $N=1$, the following elements define a basis of
$\HH^n(\Lambda_1)$ for $n \pgq 1$.
\begin{enumerate}
\item For $n=2p$ even, $n \pgq 2$:
\begin{enumerate}
\item For $-p\ppq \alpha\ppq p$: $\chi_{n,\alpha}: \ e_i\ot e_{i+2\alpha}\mapsto (-1)^{(p+\alpha)i}e_i$
for $i=0, 1$;
\item For $-p\ppq \alpha\ppq p$: $\pi_{n,\alpha}: \ e_0\ot e_{2\alpha}\mapsto
a_0\oa_0$.
\end{enumerate}
\item For $n=2p+1$ odd:
\begin{enumerate}
\item For $-p-1\ppq \gamma\ppq p$: $\varphi_{n,\gamma}:\ e_i\ot e_{i+2\gamma+1}\mapsto
(-1)^{(p+\gamma)i}a_i$ for $i=0, 1$;
\item For $-p-1\ppq\beta\ppq p$: $\psi_{n,\beta}:\ e_i\ot e_{i+2\beta+1}\mapsto
(-1)^{(p+\beta+1)i}\oa_{i+1}$ for $i=0, 1$.
\end{enumerate}
\end{enumerate}
\eprop

\bigskip

It is straightforward to give liftings of these elements and we omit
the details. Further computations enable us to give generators for
$\HH^*(\Lambda_1)$.

\bt For $m=2$ and $N=1$, $\HH^*(\L_1)$ is a finitely generated
algebra with generators:
$$\begin{array}{ll}
1, a_0\oa_0, a_1\oa_1 & \mbox{in degree 0},\\
\varphi_{1,0}, \varphi_{1,-1}, \psi_{1,0}, \psi_{1,-1} & \mbox{in degree 1},\\
\chi_{2,0}, \chi_{2,1}, \chi_{2,-1} & \mbox{in degree 2}.
\end{array}$$
\et

In order to give the structure of the Hochschild cohomology ring
modulo nilpotence, it can now be verified that $1, \chi_{2,0},
\chi_{2,1}, \chi_{2,-1}$ are the non-nilpotent generators of
$\HH^*(\L_1)$ and that $\chi_{2,0}^2 = \chi_{2,1}\chi_{2,-1}$. This
gives the following result.

\bc For $m=2$ and $N=1$,
$$\HH^*(\L_1)/\N \cong \kk[\chi_{2,0}, \chi_{2,1}, \chi_{2,-1}]/(\chi_{2,0}^2 - \chi_{2,1}\chi_{2,-1})$$
with $\chi_{2,0}, \chi_{2,1}, \chi_{2,-1}$ all in degree 2. Thus
$\HH^*(\L_1)/\N$ is a commutative finitely generated algebra of
Krull dimension 2. \ec

\bigskip

\subsection{The case $N>1$}

We start by giving a basis of cocycles for each cohomology group,
together with a lifting for each basis element. We shall use once more
the notation $\rep{n}\in\set{0,1,\ldots,m-1}$ for the representative
of $n$ modulo $m.$

\bprop For $m=2$ and $N > 1$,
$$\dim\HH^n(\lan)=
\begin{cases}
2N+1 &\mbox{if $n=0,$}\\
2N+2n &\mbox {if $n \pgq 1$ and $\car\kk\nmid N$,}\\
2N+2n+1 &\mbox {if $n \pgq 1$ and $\car\kk\mid N$.}
\end{cases}
$$
\eprop

\bprop For $m=2$, $N>1$ and for each $n \pgq 1$, we give a basis for
the cohomology group $\HH^n(\Lambda_N)$ together with a lifting for
each basis element.
\begin{enumerate}
\item For $n$ even, $n \pgq 2$:
\begin{enumerate}[$\bullet$]
\item For $-p\ppq \alpha \ppq p$:
  \begin{enumerate}[$\diamond$]
  \item $\chi_{n,\alpha}:\ e_i\ot e_{i+2\alpha }\mapsto (-1)^{\left(\frac{n}{2}-\alpha\right)i}e_i$ for $i=0, 1$,
  \item ${\mathcal L}^q{\chi}_{n,\alpha}(e_i\ot e_{i+q-2\ell+2\alpha})=\\
  \begin{cases}
  (-1)^{\left(\frac{n}{2}-\alpha\right)i}e_i\ot e_{i+q-2\ell} &\mbox{if $(q-2\ell)\alpha\pgq 0$}\\
  0&\mbox{if $(q-2\ell)\alpha<0$ and $\abs{q-2\ell}>2$}\\
  (-1)^{\left(\frac{n}{2}-\alpha\right)i} (a\oa)^{N-1}e_i\ot e_{i+2}(\oa a)^{N-1}
 & \mbox{if $\alpha<0$ and $q-2\ell=2$}\\
  (-1)^{\left(\frac{n}{2}-\alpha\right)i} (\oa a)^{N-1}e_i\ot e_{i-2}(a\oa )^{N-1}
 & \mbox{if  $\alpha>0$ and $q-2\ell=-2$}\\
  (-1)^{\left(\frac{n}{2}-\alpha\right)i} \left[\sum_{k=0}^{N-1}(a\oa)^ke_i\ot
e_{i+1}(\oa a)^{N-k-1}+\right.\\
\hspace*{1cm}\left.(-1)^{\frac{q+1}{2}}\sum_{k=0}^{N-2}(a\oa)^ka e_{\rep{i+1}}\ot e_{\rep{i+1}-1} a(\oa a)^{N-k-2}\right]  &\mbox{if  $\alpha<0$ and $q-2\ell=1$}\\
  (-1)^{\left(\frac{n}{2}-\alpha\right)i} \left[\sum_{k=0}^{N-1}(\oa a)^ke_i\ot
e_{i-1}(a\oa)^{N-k-1}+\right.\\
\hspace*{1cm}\left.(-1)^{\frac{q+1}{2}}\sum_{k=0}^{N-2}(\oa a)^k\oa e_{\rep{i-1}}\ot e_{\rep{i-1}+1} \oa(a\oa)^{N-k-2}\right]  &\mbox{if $\alpha>0$ and $q-2\ell=-1$}\\
  \end{cases}$\\
  for $i=0, 1$ and $0\ppq\ell\ppq q$.
  \end{enumerate}

\item For $-p\ppq \alpha \ppq p$:
  \begin{enumerate}[$\diamond$]
  \item $\pi_{n,\alpha}:\ e_0\ot e_{2\alpha }\mapsto (a_0\oa_0)^N,$
  \item ${\mathcal L}^q{\pi}_{n,\alpha}(e_0\ot e_{q-2\ell+2\alpha })=
  \begin{cases}
  (a\oa)^N e_0\ot e_{q-2\ell} &\mbox{if $\alpha(q-2\ell)\pgq 0$ }\\
  (a\oa)^Ne_0\ot e_{1}(\oa a)^{N-1}&\mbox{if $q-2\ell=1$ and $\alpha<0$}\\
  (a\oa)^Ne_0\ot e_{-1}(a\oa)^{N-1}&\mbox{if $q-2\ell=-1$ and $\alpha>0$}\\
  0&\mbox{otherwise}
  \end{cases}$\\
  for all $0\ppq\ell\ppq q$.
  \end{enumerate}

\item For each $1\ppq s\ppq N-1$:
  \begin{enumerate}[$\diamond$]
  \item $F_{n,0,s}:
  \begin{cases}
  e_0\ot e_0\mapsto (a_0\oa_0)^s\\e_{1}\ot e_{1}\mapsto
  (-1)^{\frac{n}{2}}(\oa_0a_0)^s,
  \end{cases}$
  \item ${\mathcal L}^q{F}_{n,0,s}:
  \begin{cases}
  e_0\ot e_{q-2\ell}\mapsto (a_0\oa_0)^se_0\ot
  e_{q-2\ell}\\e_{1}\ot e_{1+q-2\ell}\mapsto
  (-1)^{\frac{n}{2}}(\oa_0a_0)^s e_{1}\ot e_{1+q-2\ell},
  \end{cases}$\\
  for all $0\ppq\ell\ppq q$.
  \end{enumerate}

\item For each $1\ppq s\ppq N-1$:
  \begin{enumerate}[$\diamond$]
  \item $F_{n,1,s}:
  \begin{cases}
  e_1\ot e_1\mapsto (a_1\oa_1)^s\\e_0\ot e_0\mapsto
  (-1)^{\frac{n}{2}}(\oa_1a_1)^s,
  \end{cases}$
  \item ${\mathcal L}^q{F}_{n,1,s}:
  \begin{cases}
  e_1\ot e_{1+q-2\ell}\mapsto (a_1\oa_1)^se_1\ot
  e_{1+q-2\ell}\\e_0\ot e_{q-2\ell}\mapsto
  (-1)^{\frac{n}{2}}(\oa_1a_1)^s e_0\ot e_{q-2\ell},
  \end{cases}$\\
  for all $0\ppq\ell\ppq q$.
  \end{enumerate}

\item Additionally for $\car\kk\mid N$:
  \begin{enumerate}[$\diamond$]
  \item $\theta_{n}:e_1\ot e_1\mapsto (a_1\oa_1)^N,$
  \item ${\mathcal L}^q{\theta}_{n}:e_1\ot e_{1+q-2\ell}\mapsto (a_1\oa_1)^Ne_1\ot e_{1+q-2\ell}$
  for all $0\ppq\ell\ppq q$.
  \end{enumerate}
\end{enumerate}

\item For $n$ odd:
\begin{enumerate}[$\bullet$]
\item For $-p-1\ppq \gamma<0$:
  \begin{enumerate}[$\diamond$]
  \item $\varphi_{n,\gamma}:\ e_i\ot e_{i+2\gamma +1}\mapsto
  (-1)^{\left(\frac{n-1}{2}-\gamma \right)i}a_{i}(\oa a)^{N-1}$ for $i=0, 1,$
  \item ${\mathcal L}^q{\varphi}_{n,\gamma}(e_i\ot e_{i+2\gamma +1+q-2\ell})=
  \begin{cases}
  0&\mbox{if $q-2\ell>1$}\\
  -(-1)^{\left(\frac{n-1}{2}-\gamma \right)i}(a\oa)^{N-1} e_i\ot e_{i+1}a(\oa a)^{N-1}&\mbox{if $q-2\ell=1$}\\
  (-1)^{\left(\frac{n-1}{2}-\gamma \right)i}  e_i\ot
  e_{i+q-2\ell}a(\oa a)^{N-1}&\mbox{if $q-2\ell\ppq 0$}
  \end{cases}$\\
  for $i=0, 1$ and for all $0\ppq\ell\ppq q$.
  \end{enumerate}

\item For $0\ppq \gamma \ppq p:$
  \begin{enumerate}[$\diamond$]
  \item $\varphi_{n,\gamma}:\ e_i\ot e_{i+2\gamma +1}\mapsto (-1)^{\left(\frac{n-1}{2}-\gamma\right)i}a_{i}$ for $i=0, 1,$
  \item If $\gamma>0$, then ${\mathcal L}^q{\varphi}_{n,\gamma}(e_i\ot e_{i+2\gamma +1+q-2\ell})=\\
  \begin{cases}
  (-1)^{q+\left(\frac{n-1}{2}-\gamma \right)i} e_i\ot e_{i+q-2\ell}a &\mbox{if $q-2\ell\pgq 0$}\\
  -(-1)^{q+\left(\frac{n-1}{2}-\gamma \right)i}
    \left[\dsty\sum_{k=0}^{N-1}(\oa a)^ke_i\ot e_{i-1}(a\oa)^{N-k-1}a+\right.\\
    \hspace*{3cm}\left.(-1)^{\frac{q+1}{2}}\dsty\sum_{k=0}^{N-2}(\oa a)^k\oa e_{\rep{i-1}}\ot e_{\rep{i-1}+1} (\oa a)^{N-k-1}\right]
    &\mbox{if $q-2\ell=-1$}\\
  (-1)^{q+\left(\frac{n-1}{2}-\gamma \right)i}(\oa a)^{N-1}e_i\ot e_{i-2}(a\oa)^{N-1}a &\mbox{if $q-2\ell=-2$}\\
  0&\mbox{if $q-2\ell<-2$}
  \end{cases}$\\
  for $i=0, 1$ and for all $0\ppq\ell\ppq q$.

  \item If $\gamma=0$, then ${\mathcal L}^q{\varphi}_{n,0}(e_i\ot e_{i+1+q-2\ell})=\\
  \begin{cases}
  (-1)^{q+\frac{n-1}{2}i} e_i\ot e_{i+q-2\ell}a   &\mbox{if $q-2\ell=q$}\\
  (-1)^{q+\frac{n-1}{2}i} \left[e_i\ot e_{i+q-2\ell}a -(-1)^q(N-1)e_i\ot e_{i+q-2\ell+2}\oa(a\oa)^{N-1} \right] &\mbox{if $0\ppq q-2\ell<q$} \\
  (-1)^{q+\frac{n-1}{2}i} N e_i\ot e_{i+q-2\ell}a(\oa a)^{N-1}   &\mbox{if $q-2\ell<-1$} \\
  \!\!\!\begin{array}{l}
  (-1)^{1+\frac{n-1}{2}i} \left\{\dsty\sum_{k=1}^{N-1}\sum_{v=0}^{k-1}
    \left[(a\oa)^ve_i\ot e_{i+1}\oa(a\oa)^{N-v-1}+\right.\right.\\
    \hspace*{3.5cm}\left.(-1)^{\frac{q+1}{2}}(a\oa)^vae_{\rep{i+1}}\ot e_{\rep{i+1}-1}(a\oa)^{N-v-1}+\right.\\
    \hspace*{4cm}\left.(-1)^{\frac{q+1}{2}}(\oa a)^v\oa e_{\rep{i-1}}\ot e_{\rep{i-1}+1}(\oa a)^{N-v-1}+\right.\\
    \hspace*{4.5cm}\left.(\oa a)^ve_i\ot e_{i-1}a(\oa a)^{N-v-1}\right]\\
  \left.\hspace*{2.5cm}+\dsty\sum_{k=0}^{N-1}(\oa a)^ke_i\ot e_{i-1}
  a(\oa a)^{N-k-1}\right\}
  \end{array} &\mbox{if $q-2\ell=-1$}
  \end{cases}$\\
  for $i=0, 1$ and for all $0\ppq\ell\ppq q$.
  \end{enumerate}

\item For $0\ppq\beta\ppq p$:
  \begin{enumerate}[$\diamond$]
  \item $\psi_{n,\beta}:e_i\ot e_{i+2\beta +1}\mapsto (-1)^{\left(\frac{n+1}{2}-\beta\right)i}\oa(a\oa)^{N-1},$
  \item ${\mathcal L}^q{\psi}_{n,\beta}(e_i\ot e_{i+2\beta +1+q-2\ell})=
  \begin{cases}
  0&\mbox{if $q-2\ell<-1$}\\
  (-1)^{\left(\frac{n+1}{2}-\beta\right)i}(\oa a)^{N-1}e_i\ot e_{i-1}\oa(a\oa)^{N-1}&\mbox{if $q-2\ell=-1$}\\
  (-1)^{\left(\frac{n+1}{2}-\beta\right)i} e_i\ot
  e_{i+q-2\ell}\oa(a\oa)^{N-1}&\mbox{if $q-2\ell>-1$}
  \end{cases}$\\
  for $i=0, 1$ and for all $0\ppq\ell\ppq q$.
  \end{enumerate}

\item For $-p-1\ppq \beta < 0$:
  \begin{enumerate}[$\diamond$]
  \item $\psi_{n,\beta}:e_i\ot e_{i+2\beta +1}\mapsto (-1)^{\left(\frac{n+1}{2}-\beta\right)i}\oa,$
  \item If $\beta < -1,$ then ${\mathcal L}^q{\psi}_{n,\beta}(e_i\ot e_{i+2\beta +1+q-2\ell})=\\
  \begin{cases}
  (-1)^{\left(\frac{n+1}{2}-\beta\right)i}  e_i\ot e_{i+q-2\ell}\oa     &\mbox{if $q-2\ell<1$}\\
  (-1)^{\left(\frac{n+1}{2}-\beta\right)i}
     \left[\dsty\sum_{k=0}^{N-1}(a\oa)^k e_i\ot e_{i+1}\oa (a\oa)^{N-k-1}+\right.\\
     \hspace*{3cm}\left.\dsty (-1)^{\frac{q+1}{2}}\sum_{k=0}^{N-2}(a\oa)^kae_{\rep{i+1}}\ot e_{\rep{i+1}-1}(a\oa)^{N-k-1}\right]
     &\mbox{if $q-2\ell=1$}\\
  (-1)^{\left(\frac{n+1}{2}-\beta\right)i} (a\oa)^{N-1}e_i\ot e_{i+2}(\oa a)^{N-1}\oa      &\mbox{if $q-2\ell=2$}\\
   0      &\mbox{if $q-2\ell>2$}
  \end{cases}$\\
  for $i=0, 1$ and for all $0\ppq\ell\ppq q$.

  \item If $\beta=-1$, then ${\mathcal L}^q{\psi}_{n,-1}(e_i\ot e_{i-1+q-2\ell})=\\
  \begin{cases}
  (-1)^{\frac{n-1}{2}i} e_i\ot e_{i+q-2\ell}\oa   &\mbox{if $q-2\ell=-q$}\\
  (-1)^{\frac{n-1}{2}i} \left[e_i\ot e_{i+q-2\ell}\oa -(-1)^q(N-1)e_i\ot e_{i+q-2\ell-2}a(\oa a)^{N-1} \right] &\mbox{if $-q< q-2\ell\ppq 0$} \\
  (-1)^{\frac{n-1}{2}i} N e_i\ot e_{i+q-2\ell}\oa(a \oa)^{N-1}   &\mbox{if $q-2\ell>1$} \\
  \begin{array}{l}
  \!\!\!(-1)^{\frac{n-1}{2}i} \left\{\dsty\sum_{k=1}^{N-1}\sum_{v=0}^{k-1}
    \left[(\oa a)^ve_i\ot e_{i-1}a(\oa a)^{N-v-1}+\right.\right.\\
    \hspace*{3cm}\left.(-1)^{\frac{q+1}{2}}(\oa a)^v\oa e_{\rep{i-1}}\ot e_{\rep{i-1}+1}(\oa a)^{N-v-1}+\right.\\
    \left.\hspace*{3.5cm}(-1)^{\frac{q+1}{2}}(a \oa)^va e_{\rep{i+1}}\ot e_{\rep{i+1}-1}(a \oa)^{N-v-1}+\right.\\
    \hspace*{4cm}\left.(a \oa)^ve_i\ot e_{i+1}\oa(a \oa)^{N-v-1}\right]\\
  \left.\hspace*{2cm}+\dsty\sum_{k=0}^{N-1}(a \oa)^ke_i\ot e_{i+1}
  \oa(a \oa)^{N-k-1}\right\}
  \end{array} &\mbox{if $q-2\ell=1$}
  \end{cases}$\\
  for $i=0, 1$ and for all $0\ppq\ell\ppq q$.
  \end{enumerate}

\item For each $j = 0, 1$ and each $1\ppq s\ppq N-1$:
  \begin{enumerate}[$\diamond$]
  \item $E_{n,j,s}:e_j\ot e_{j-1}\mapsto (\oa a)^s\oa_j$,
  \item ${\mathcal L}^q{E}_{n,j,s}:\\
  \begin{cases}
  e_j\ot e_{j-1}\mapsto  (\oa a)^s e_j\ot e_j \oa &\mbox{if $q$ is even}\\
  e_j\ot e_j\mapsto \dsty\sum_{k=0}^{N-s}\sum_{v=k}^{N-s}
    \left[(\oa a)^{v+s} e_j\ot e_{j-1}a(\oa a)^{N-v-1}+\right.\\
    \hspace*{3cm}\left.(-1)^{\frac{q+1}{2}}(\oa a)^{v+s}ae_{j-1}\ot e_{j}(\oa a)^{N-v-1}\right] &\mbox{if $q$ is odd}\\
  \!\!\!\begin{array}{ll}e_{\rep{j-1}}\ot e_{\rep{j-1}}\mapsto
  (-1)^{\frac{n-1}{2}}\dsty\sum_{k=0}^{N-s-1}\sum_{v=k}^{N-s}&\left[(-1)^{\frac{q+1}{2}}(a
  \oa)^{v+s}a e_j\ot e_{j+1}(a \oa)^{N-v-1}\right.\\&\left.+(a
  \oa)^{v+s+1}e_{\rep{j+1}}\ot e_{\rep{j+1}-1} \oa(a \oa)^{N-v-2}\right]\end{array}
  &\mbox{if $q$ is odd}
  \end{cases}$\\
  for all $0\ppq\ell\ppq q$.
  \end{enumerate}

\item Additionally for $\car\kk\mid N$:
  \begin{enumerate}[$\diamond$]
  \item $\omega_{n}:e_0\ot e_{1}\mapsto a_0,$
  \item ${\mathcal L}^q\omega_{n}:\\
  \begin{cases}
  e_0\ot e_0 \mapsto \dsty\sum_{k=1}^{N-1}\sum_{v=0}^{k-1}
    \left[(-1)^{\frac{q-1}{2}}(a\oa)^vae_{1}\ot e_0(a\oa)^{N-v-1}-\right.\\
    \hspace*{3.5cm}\left.(a\oa)^ve_0\ot
      e_{1}\oa(a\oa)^{N-v-1}\right]\mbox{ if $q$ is odd}\\
  e_{1}\ot e_{1}\mapsto (-1)^{\frac{n+1}{2}}\dsty\sum_{k=1}^{N-1}\sum_{v=0}^{k-1}
    \left[(\oa a)^{v+1}e_{1}\ot e_0 a(\oa a)^{N-v-2}+\right.\\
    \hspace*{3.5cm}\left.(-1)^{\frac{q+1}{2}}(\oa a)^v\oa e_0\ot
      e_{1}(\oa a)^{N-v-1}\right]  \mbox{ if $q$ is odd}\\
  e_{\rep{-v}}\ot e_{\rep{-v}+2v}\mapsto (-1)^{\frac{n-1}{2}v}
    \left[(-1)^{v+\frac{q+1}{2}}ae_{\rep{-v+1}}\ot e_{\rep{-v+1}+2v-1}+\right.\\
    \hspace*{3.5cm}\left.e_{\rep{-v}}\ot e_{\rep{-v}+2v+1}\oa
      (a\oa)^{N-1}\right] \mbox{ for $1\ppq v\ppq \frac{q+1}{2}$ and $q$ odd}\\
  e_{\rep{-v}}\ot e_{\rep{-v}+2v+1}\mapsto (-1)^{\frac{n-1}{2}v}
  \left[e_{\rep{-v}}\ot e_{\rep{-v}+2v}a-(-1)^{v+\frac{q}{2}}\oa(a\oa)^{N-1}e_{\rep{-v-1}}\ot e_{\rep{-v-1}+2v+2}\right]
  \\ \hspace*{8cm}\mbox{for $0\ppq v\ppq \frac{q}{2}$ and $q$ even}
  \end{cases}\\$
  for all $0\ppq\ell\ppq q$.
  \end{enumerate}
\end{enumerate}
\end{enumerate}
\eprop

\bt\label{gens:m=2} For $m=2$ and $N > 1$, $\HH^*(\L_N)$ is a
finitely generated algebra, with generators:
$$\begin{array}{ll}
1, (a_0\oa_0)^N, (a_1\oa_1)^N, (a_0\oa_0 + \oa_0 a_0), (a_1\oa_1 +
\oa_1 a_1) & \mbox{in degree 0},\\
\varphi_{1,0}, \varphi_{1,-1}, \psi_{1,0}, \psi_{1,-1}, E_{1,0,1} &
\mbox{in degree 1},\\
\mbox{$\omega_{1}$ if $\car\kk\mid N$} & \mbox{in degree 1},\\
\chi_{2,0}, \chi_{2,1}, \chi_{2,-1}, F_{2,0,1} & \mbox{in degree 2}.
\end{array}$$
\et

It is easy to verify that the non-nilpotent generators of
$\HH^*(\L_N)$ are $1, \chi_{2,0}, \chi_{2,1}$ and $\chi_{2,-1}$, and
that $\chi_{2,1}\chi_{2,-1} = 0$. Thus we have the following
theorem.

\bc For $m=2$ and $N > 1$,
$$\HH^*(\L_N)/\N \cong \kk[\chi_{2,0}, \chi_{2,1}, \chi_{2,-1}]/(\chi_{2,1}\chi_{2,-1}).$$
Hence $\HH^*(\L_N)/\N$ is a commutative finitely generated algebra
of Krull dimension 2. \ec

\bigskip

\section{The Hochschild cohomology ring $\HH^*(\L_N)$ for $m=1$}

We start by giving the centre of the algebra, since $\HH^0(\L_N) =
Z(\L_N)$.

\bt\label{thm:m1centre} Suppose that $m=1$. Then the dimension of
$\HH^0(\lan)$ is $N+3$, and a $\kk$-basis of $\HH^0(\lan)$ is given
by
$$\set{1,(a\oa)^N, [(a\oa)^s+(\oa a)^s] \mbox{ for $1\ppq s\ppq
N-1$}, a(\oa a)^{N-1},\oa(a\oa)^{N-1}}.$$ Thus $\HH^0(\lan)$ is
generated \textit{as an algebra} by
$$\begin{array}{ll}
\set{1,a, \oa} & \mbox{if $N=1$;} \\
\set{1,(a\oa)^N, [a\oa+\oa a], a(\oa a)^{N-1},\oa(a\oa)^{N-1}} &
\mbox{if $N>1$.} \end{array}$$\et

\subsection{The case $N=1$}

The algebra $\L_1$ is a commutative Koszul algebra. The quiver has a
single vertex $1$ and two loops $a$ and $\oa$ such that $a^2=0,$
$\oa^2=0$, and $a\oa=\oa a.$ The structure of the Hochschild
cohomology ring for $\L_1$ with $m=1$ was determined in \cite{BGMS};
we state the result here for completeness. Note that we have
$\HH^0(\L_1) = \L_1$ in this case.

\bt \cite{BGMS} For $m=1$ and $N=1$, we have $$\HH^*(\L_1) \cong
\begin{cases}
\L_1\langle u_0, u_1, x_0, x_1\rangle/\langle u_0^2, u_1^2, au_0, \oa u_1,
ax_0, \oa x_1\rangle & \mbox{if $\car\kk\neq 2$,}\\
\L_1[y_0,y_1] & \mbox{if $\car\kk=2$,}
\end{cases}
$$ where
$x_0, x_1$ are in degree 2 and $u_0, u_1, y_0, y_1$ are in degree 1.
\et

We remark that we could also have used the K\"unneth formula
\cite[Theorem 7.4 p297]{ML} together with $\HH^*(A)$, where
$A=\kk[x]/(x^2)$. The cohomology of $A$ has been studied in several places
in the literature including in \cite{C,L}.

We can now give the structure of the Hochschild cohomology ring of $\L_1$
modulo nilpotence and its Krull dimension.

\bc For $m=1$ and $N=1$,
$$\HH^*(\la)/\mathcal{N}\cong
\begin{cases}
\kk[x_0, x_1] & \mbox{if $\car\kk\neq 2$}\\
\kk[y_0, y_1] & \mbox{if $\car\kk = 2$,}
\end{cases}
$$
where $x_0, x_1$ are in degree 2 and $y_0, y_1$ are in degree 1.
Hence $\HH^*(\la)/\N$ is a commutative finitely generated algebra of
Krull dimension $2$. \ec

\bigskip

\subsection{The case $N>1$ and $\car\kk\neq 2$}

For the case $N > 1$ and $m=1$ there are two subcases to consider,
depending on whether the characteristic of the field is or is not
equal to 2. Firstly we suppose that $\car\kk \neq 2$.

\bprop For $m=1, N>1$ and $\car\kk \neq 2$, the dimensions of the
Hochschild cohomology groups are given by $\dim\HH^0(\lan) = N+3$
and, for $n \pgq 1$,
$$\dim\HH^n(\lan) = \left \{ \begin{array}{ll}
N+4n+3& \mbox{if $\car\kk=2$, }\\
N+n+2 & \mbox{if $\car\kk\neq2$, $\car\kk\nmid N$,}\\
N+n+2 & \mbox{if $\car\kk\neq2$, $\car\kk\mid N$, $n\md{1,2}$},\\
N+n+3 & \mbox{if $\car\kk\neq2$, $\car\kk\mid N$, $n\md{0,3}$.}
\end{array}\right.
$$
\eprop

We now describe explicitly the elements of the Hochschild cohomology
groups $\HH^n(\lan)$ for $n \pgq 1$, in order to give a set of
generators for $\HH^*(\lan)$ as a finitely generated algebra. Recall
that $\HH^0(\lan) = Z(\lan)$ and was described in Theorem
\ref{thm:m1centre}. For each cocycle $f\in\Hom(P^{n},\lan)$, we
continue our practice of writing only the image of each generator
$e_i\ot e_{i+n-2r}$ in $P^n$ where that image is non-zero. However,
since we have a single vertex, we write 1 for $e_0$, and we use the
notation $1 \ot_r 1$ instead of $e_0\ot e_{n-2r}$ for the generator
of the $r$-th summand $\lan e_0\ot e_{n-2r}\lan$ of $P^n$. (Recall
that $r = 0, \ldots , n$.) This notation ensures that, in this
single idempotent case $m=1$, we continue to avoid any confusion as
to which summand of $P^n$ we are referring.

\bprop Suppose that $m=1, N>1$ and $\car\kk \neq 2$. The following
elements define a basis of $\HH^n(\lan)$ for $n \pgq 1$.
\begin{enumerate}
\item For all $n \pgq 1$ and $0\ppq r\ppq n$ with $r$ even: $\chi_{n,r}: 1\ot_r 1\mapsto
1$.
\item For $n$ even, $n \pgq 2$:
\begin{enumerate}
\item $\psi_{n,0}:1\ot_0 1\mapsto \oa(a\oa)^{N-1}$;
\item $\varphi_{n,n}:1\ot_n 1\mapsto a(\oa a)^{N-1}$;
\item For $1\ppq r\ppq n$ and $r$ odd: $\pi_{n,r}:1\ot_r 1\mapsto (a\oa)^N$;
\item For each $s$ such that $1\ppq s\ppq N-1$: $F_{n,s}: 1\ot_{n/2} 1
\mapsto (a\oa)^s+ (-1)^\frac{n}{2}(\oa a)^s$;
\item Additionally, for $\car\kk\mid N$, if $n \equiv 0 \pmod 4$:
$\theta_{n}:1\ot_{n/2} 1\mapsto (a\oa)^N$.
\end{enumerate}
\item For $n$ odd:
\begin{enumerate}
\item For $0\ppq r\ppq \frac{n-1}{2}$ and $r$ even: $\varphi_{n,r}:1\ot_r 1\mapsto a$;
\item For $\frac{n+1}{2}\ppq r\ppq n-1$ and $r$ even: $\varphi_{n,r}:1\ot_r 1\mapsto
(a\oa)^{N-1}a$;
\item For $1\ppq r\ppq \frac{n-1}{2}$ and $r$ odd: $\psi_{n,r}:\ 1\ot_r 1\mapsto
(\oa a)^{N-1}\oa$;
\item For $\frac{n+1}{2}\ppq r\ppq n$ and $r$ odd: $\psi_{n,r}:\ 1\ot_r 1\mapsto
\oa$;
\item $\pi_{n,+}:1\ot_0 1\mapsto (a\oa)^N$;
\item $\pi_{n,-}:1\ot_n 1\mapsto (a\oa)^N$;
\item For each $s$ such that $1\ppq s\ppq N-1$: $E_{n,s}:1\ot_{(n-1)/2} 1\mapsto (a\oa)^sa$.
\item Additionally, for $\car\kk\mid N$, if $n \equiv 3 \pmod 4$ and $r=\frac{n+1}{2}$:
$\psi_{n,r}:1\ot_r 1\mapsto \oa $.
\end{enumerate}
\end{enumerate}
\eprop

We omit the liftings for this case, since the calculations are
straightforward, and give generators of the Hochschild cohomology
ring.

\bt For $m=1, N>1$, $\car\kk \neq 2$, $\HH^*(\L_N)$ is a finitely
generated algebra with generators: $$\begin{array}{ll} 1, a\oa+\oa
a, a(\oa a)^{N-1}, \oa(a\oa)^{N-1}, (a\oa)^N & \mbox{in
degree 0},\\
\varphi_{1,0}, \psi_{1,1} & \mbox{in degree 1},\\
\chi_{2,0}, \chi_{2,2}, \pi_{2,1}, F_{2,1}, \varphi_{2,2},
\psi_{2,2} & \mbox{in degree 2},\\
\varphi_{3,2}, \psi_{3,1} & \mbox{in degree 3},\\
\mbox{$\psi_{3,2}$ if $\car\kk\mid N$} & \mbox{in degree 3},\\
\chi_{4,2}, F_{4,1}& \mbox{in degree 4.}
\end{array}$$
\et

The non-nilpotent generators of $\HH^*(\lan)$ are precisely those
elements given in the above theorem whose image does not lie in
$\rrad$, that is, the elements $1, \chi_{2,0}, \chi_{2,2}$ and
$\chi_{4,2}$. This gives the next result.

\bc For $m=1, N > 1$ and $\car\kk \neq 2$,
$$\HH^*(\lan)/\mathcal{N}\cong
\kk[\chi_{2,0}, \chi_{2,2}, \chi_{4,2}]/(\chi_{2,0}\chi_{2,2}).$$
Hence $\HH^*(\lan)/\mathcal{N}$ is a commutative finitely generated
algebra of Krull dimension 2. \ec

\bigskip

\subsection{The case $N>1$ and $\car\kk = 2$}

In this final case, we continue to keep the notation that $1 \ot_r 1
= e_0\ot e_{n-2r}$ for $r = 0, 1, \ldots , n$.

\bprop For $m=1$, $N > 1$ and $\car\kk = 2$, then $\dim\HH^{n}(\L_N)
= 4n+N+3$ for all $n \pgq 0$. \eprop

\bprop Suppose that $m=1$, $N > 1$ and $\car\kk = 2$. For $n \pgq
0$, the following elements define a basis of $\HH^n(\L_N)$.
\begin{enumerate}
\item For all $n\pgq 0$:
\begin{enumerate}
\item For $r = 0, \ldots , n$: $\chi_{n,r}:\ 1\ot_r 1\mapsto 1$;
\item For $r = 0, \ldots , n$: $\pi_{n,r}:\ 1\ot_r 1\mapsto (a\oa)^N$.
\end{enumerate}
\item For $n$ even, $n \pgq 2$:
\begin{enumerate}
\item For $r = 0, 1, \ldots \frac{n-2}{2}$: $\varphi_{n,r}:\ 1\ot_r 1\mapsto a$;
\item For $r = \frac{n}{2}, \ldots , n-1, n$: $\varphi_{n,r}:\ 1\ot_r 1\mapsto
(a\oa)^{N-1}a$;
\item For $r= 0, 1, \ldots , \frac{n}{2}$:
$\psi_{n,r}:\ 1\ot_r 1\mapsto (\oa a)^{N-1}\oa$;
\item For $r = \frac{n+2}{2}, \ldots , n-1, n$: $\psi_{n,r}:\ 1\ot_r 1\mapsto
\oa$;
\item For $j = 1, \ldots , N-1$: $F_{n,j}:\ 1\ot_{n/2}1\mapsto
(a\oa)^j + (\oa a)^j$.
\end{enumerate}
\item For $n$ odd:
\begin{enumerate}
\item For $r = 0, 1, \ldots \frac{n-1}{2}$: $\varphi_{n,r}:\ 1\ot_r 1\mapsto
a$;
\item For $r = \frac{n+1}{2}, \ldots , n-1, n$: $\varphi_{n,r}:\ 1\ot_r 1\mapsto
(a\oa)^{N-1}a$;
\item For $r = 0, 1, \ldots , \frac{n-1}{2}$: $\psi_{n,r}:\ 1\ot_r 1\mapsto (\oa
a)^{N-1}\oa$;
\item For $r = \frac{n+1}{2}, \ldots , n-1, n$: $\psi_{n,r}:\ 1\ot_r 1\mapsto
\oa$;
\item For $j = 1, \ldots, N-1$: $E_{n,j}:\ 1\ot_{\frac{n+1}{2}}1\mapsto \oa(a\oa)^j$.
\end{enumerate}
\end{enumerate}
\eprop

Again we omit the details of the liftings, and state directly the
generators of $\HH^*(\L_N)$.

\bt For $m=1, N > 1$ and $\car\kk = 2$, the Hochschild cohomology
ring $\HH^*(\L_N)$ is finitely generated as a $\kk$-algebra with
generators
$$\begin{array}{ll}
1, a\oa + \oa a, (a\oa)^{N-1}a, (\oa a)^{N-1}\oa, (a\oa)^N & \mbox
{in degree 0},\\
\chi_{1,0}, \chi_{1,1}, \varphi_{1,0}, \psi_{1,1} & \mbox {in degree 1},\\
\chi_{2,1}, \psi_{2,2} & \mbox {in degree 2}.
\end{array}$$
\et

\bc For $m=1, N > 1$ and $\car\kk = 2$,
$$\HH^*(\L_N)/\N \cong \kk[\chi_{1,0}, \chi_{1,1}, \chi_{2,1}]/(\chi_{1,0}\chi_{1,1})^2.$$
Thus $\HH^*(\L_N)/\N$ is a commutative ring of Krull dimension 2.
\ec

\bigskip

\section{Summary}

In conclusion we have the following theorem.

\bt For all $N \pgq 1$ and $m \pgq 1$, $\HH^*(\L_N)$ is a finitely
generated $\kk$-algebra. Moreover $\HH^*(\L_N)/\N$ is a finitely
generated $\kk$-algebra of Krull dimension 2. \et

In particular the theorem holds for the algebras $\L_1$ which occur
in the representation theory of the Drinfeld doubles of the
generalised Taft algebras (\cite{EGST}), and for the algebras of
Farnsteiner and Skowro\'nski (\cite{FS}). Moreover the conjecture of
\cite{SS} has been proved for this class of algebras.

\bigskip

We conclude the paper by recalling the following simultaneous conditions which
were considered in \cite{EHSST}:\\
{\bf (Fg1)}\ There is a commutative Noetherian graded subalgebra $H$
of $\HH^*(\L_1)$ such that $H^0 = \HH^0(\L_1)$.\\
{\bf (Fg2)}\ $\Ext^*_{\L_1}(\L_1/\rrad, \L_1/\rrad)$ is a finitely
generated $H$-module.

In the case $N=1$, it is immediate from the above results that the
algebras $\L_1$ satisfy the condition (Fg1) with
$H=\HH^{\ev}(\L_1)$. We know that $\L_1$ is a Koszul algebra with
$\Ext^*_{\L_1}(\L_1/\rrad, \L_1/\rrad)$ being the Koszul dual of
$\L_1$. Hence $\Ext^*_{\L_1}(\L_1/\rrad, \L_1/\rrad) = \kk{\mathcal
Q}/T$ where ${\mathcal Q}$ is the quiver of $\L_1$ and $T$ is the
ideal of $\kk{\mathcal Q}$ generated by $a_i\oa_i +
\oa_{i-1}a_{i-1}$ for all $i =0, 1, \ldots , m-1$. It may be
verified that the condition (Fg2) also holds with
$H=\HH^{\ev}(\L_1)$, the full details being left to the reader.

The conjecture of \cite{SS} concerning the finite generation of the
Hochschild cohomology ring modulo nilpotence was motivated by the
work on support varieties. The fact that $\L_1$, which is
self-injective, and $H = \HH^{\ev}(\L_1)$ together satisfy (Fg1) and
(Fg2) means that the support varieties for finitely generated
$\L_1$-modules satisfy all the properties of \cite{SS} and the
subsequent work of \cite{EHSST} on self-injective algebras.

\section*{Acknowledgements}

The authors thank Karin Erdmann for all her helpful comments and
valuable suggestions on this paper.


\bigskip

\begin{thebibliography}{99}

\bibitem{BGMS} Buchweitz, R.-O., Green, E.L., Madsen, D. and
Solberg, \O., {\em Finite Hochschild cohomology without finite
global dimension}, Math. Res. Lett. {\bf 12} (2005), 805-816.

\bibitem{BGSS} Buchweitz, R.-O., Green, E.L., Snashall, N. and
Solberg, \O., {\em Multiplicative structures for Koszul algebras},
Quart. J. Math. (to appear).

\bibitem{BK} Butler, M.C.R and King, A.D., {\em Minimal resolutions of
algebras}, J. Algebra {\bf 212} (1999), 323-362.

\bibitem{C} Cibils, C., {\em Hochschild cohomology algebra of
radical square zero algebras,} Algebras and modules, II (Geiranger, 1996),  93-101,
CMS Conf. Proc., {\bf 24}, Amer. Math. Soc., Providence, RI, 1998.

\bibitem{CK} Chin, W. and Krop, L., {\em Representation theory of liftings of
quantum planes,} arXiv:0712.1078.

\bibitem{EGST} Erdmann, K., Green, E.L., Snashall, N. and Taillefer,
R., {\em Representation theory of the Drinfeld doubles of a family
of Hopf algebras}, J. Pure Appl. Algebra {\bf 204} (2006), 413-454.

\bibitem{EHSST} Erdmann, K., Holloway, M., Snashall, N., Solberg, \O. and Taillefer,
R., {\em Support varieties for selfinjective algbras}, K-Theory {\bf
33} (2004), 67-87.

\bibitem{FS} Farnsteiner, R. and Skowro\'nski, A., {\em Classification of
restricted Lie algebras with tame principal block}, J. Reine Angew.
Math. {\bf 546} (2002), 1-45.

\bibitem{FS2} Farnsteiner, R. and Skowro\'nski, A., {\em The tame
infinitesimal groups of odd characteristic}, Adv. Math. {\bf 205}
(2006), 229-274.

\bibitem {GHMS} Green, E.L., Hartman, G., Marcos E.N. and Solberg, \O.,
{\em Resolutions over Koszul algebras}, Archiv der Math. {\bf 85}
(2005), 118-127.

\bibitem{GSn} Green, E.L. and Snashall, N., {\em Projective bimodule resolutions
of an algebra and vanishing of the second Hochschild cohomology group},
Forum Math. {\bf 16} (2004), 17-36.

\bibitem{GSS} Green, E.L., Snashall, N. and Solberg, \O., {\em The Hochschild
cohomology ring modulo nilpotence of a monomial algebra}, J. Algebra
Appl. {\bf 5} (2006), 153-192.

\bibitem{GSZ} Green, E.L., Solberg, \O. and Zacharia, D., {\em Minimal
projective resolutions}, Trans. Amer. Math. Soc. {\bf 353} (2001),
2915-2939.

\bibitem{H} Happel, D., {\em Hochschild cohomology of finite-dimensional
algebras}, Springer Lecture Notes in Mathematics {\bf 1404} (1989),
108-126.

\bibitem{L} Locateli, A.C. {\em Hochschild cohomology of truncated
quiver algebras}, Comm. Algebra {\bf 27} (1999), 645-664.

\bibitem{ML} MacLane, S., {\em Homology,} Classics in Mathematics, Springer-Verlag, 1995.

\bibitem{Patra} Patra, M.K., {\em On the structure of nonsemisimple
Hopf algebras}, J. Phys. A Math. Gen. {\bf 32} (1999), 159-166.

\bibitem{SS} Snashall, N. and Solberg, \O., {\em Support varieties
and Hochschild cohomology rings}, Proc. London Math. Soc. {\bf 88}
(2004), 705-732.

\bibitem{Suter} Suter, R., {\em Modules for
$\mathfrak{U}_q(\mathfrak{sl}_2)$}, Comm. Math. Phys. {\bf 163}
(1994), 359-393.

\bibitem{Xiao} Xiao, J., {\em Finite-dimensional representations of
$U_t(sl(2))$ at roots of unity}, Can. J. Math. {\bf 49} (1997),
772-787.

\end{thebibliography}
\end{document}